\documentclass[a4paper,12pt]{amsart}

\usepackage{geometry}
\usepackage{hyphenat} 
\usepackage[utf8]{inputenc}
\usepackage[T1]{fontenc}
\usepackage[all,cmtip,2cell]{xy}
\UseTwocells
\xyoption{2cell}
\usepackage{amsfonts}
\usepackage{amsthm}
\usepackage{amsmath}
\usepackage{amssymb}
\usepackage{bbm} 
\usepackage{enumerate}
\usepackage{verbatim} 
\usepackage{mathtools}
\usepackage[usenames,dvipsnames]{xcolor}
\usepackage{url}
\usepackage{chngcntr} 
\usepackage{pdflscape} 
\usepackage[titletoc]{appendix}
\usepackage{tikz-cd}
\usepackage{etoolbox}
\makeatletter
\appto{\appendices}{\def\Hy@chapapp{Appendix}}
\makeatother

\usepackage[stretch=10]{microtype}

\usepackage[pdftex]{hyperref} 
\usepackage{bookmark}

\renewcommand{\epsilon}{\varepsilon}

\theoremstyle{definition} \newtheorem{defn}[equation]{Definition}
\theoremstyle{plain} \newtheorem{teo}[equation]{Theorem}
\theoremstyle{plain} \newtheorem{lema}[equation]{Lemma}
\theoremstyle{plain} \newtheorem{prop}[equation]{Proposition}
\theoremstyle{plain} \newtheorem{corolario}[equation]{Corollary}
\theoremstyle{remark} \newtheorem{obs}[equation]{Remark}
\theoremstyle{remark} 
\theoremstyle{definition} 
\theoremstyle{remark} \newtheorem{ej}[equation]{Example}
\theoremstyle{remark}

\DeclareUnicodeCharacter{00A0}{ }

\DeclareMathAlphabet{\mathbbe}{U}{bbold}{m}{n}
\def\DDelta{{\mathbbe{\Delta}}}

\newcommand{\bq}{\begin{quote} \it}
\newcommand{\eq}{\it\end{quote}}

\newcommand{\D}{\mathcal{D}}
\renewcommand{\1}{\ensuremath{\mathbbm{1}}}
\newcommand{\V}{\mathcal{V}}

\newcommand{\W}{\mathcal{W}}

\renewcommand{\S}{\mathbb{S}}
\newcommand{\bbS}{\mathbb{S}}

\newcommand{\N}{\mathbb{N}}

\newcommand{\be}{\begin{enumerate}}
\newcommand{\ee}{\end{enumerate}}
\newcommand{\bi}{\begin{itemize}}
\newcommand{\ei}{\end{itemize}}
\newcommand{\bprop}{\begin{prop}}
\newcommand{\eprop}{\end{prop}}
\newcommand{\bej}{\begin{ej}}
\newcommand{\eej}{\end{ej}}
\newcommand{\bdefn}{\begin{defn}}
\newcommand{\edefn}{\end{defn}}
\newcommand{\bprf}{\begin{proof}}
\newcommand{\eprf}{\end{proof}}
\newcommand{\bobs}{\begin{obs}}
\newcommand{\eobs}{\end{obs}}
\newcommand{\bcor}{\begin{corolario}}
\newcommand{\ecor}{\end{corolario}}
\newcommand{\bteo}{\begin{teo}}
\newcommand{\eteo}{\end{teo}}

\newcommand{\Fun}{\mathop{\rm Fun}}

\newcommand{\id}{\mathrm{id}}

\newcommand{\Set}{\ensuremath{\mathbf{Set}}}

\newcommand{\Ab}{\ensuremath{\mathbf{Ab}}}
\newcommand{\Mon}{\ensuremath{\mathbf{Mon}}}
\newcommand{\CMon}{\ensuremath{\mathbf{CMon}}}
\newcommand{\CoComon}{\ensuremath{\mathbf{CoComon}}}
\newcommand{\Comon}{\ensuremath{\mathbf{Comon}}}

\newcommand{\biHopf}{\ensuremath{\mathbf{BiHopf}}}
\newcommand{\Mod}{\mbox{-}\mathbf{Mod}}

\newcommand{\RMod}{\ensuremath{R}\mbox{-}\mathbf{Mod}}

\newcommand{\RCAlg}{\ensuremath{R}\mbox{-}\ensuremath{\mathbf{CAlg}}}

\newcommand{\SMod}{\ensuremath{\mathbb{S}}\mbox{-}\ensuremath{\textbf{Mod}}}
\newcommand{\kMod}{\ensuremath{\mathbf{k}}\mbox{-}\ensuremath{\textbf{Mod}}}
\newcommand{\kmod}{\ensuremath{\mathbf{k}}\mbox{-}\ensuremath{\textbf{Mod}}}
\newcommand{\kdgm}{\ensuremath{\mathbf{k}}\mbox{-}\ensuremath{\mathbf{dgm}}}

\newcommand{\kcoco}{E_*\mbox{-}\ensuremath{\mathbf{CoCoalg}}}
\newcommand{\kkcoco}{k\mbox{-}\ensuremath{\mathbf{CoCoalg}}}
\newcommand{\rrcoco}{R'\mbox{-}\ensuremath{\mathbf{CoCoalg}}}
\newcommand{\rcoco}{R\mbox{-}\ensuremath{\mathbf{CoCoalg}}}

\newcommand{\NGrRing}{\ensuremath{\mathbf{Gr}\mathbf{Ring}}}
\newcommand{\NGrAb}{\ensuremath{\mathbf{Gr}\mathbf{Ab}}}
\newcommand{\Gr}{\ensuremath{\mathbf{Gr}}}
\newcommand{\Ring}{\ensuremath{\mathbf{Ring}}}

\newcommand{\Cat}{\ensuremath{\mathbf{Cat}}}
\newcommand{\Top}{\ensuremath{\mathbf{Top}}}

\newcommand{\diag}{\operatorname{diag}}

\newcommand{\Nat}{\mathrm{Nat}}

\newcommand{\op}{\mathrm{op}}

\newcommand{\sip}{\Sigma^\infty_+}
\newcommand{\aug}{\mathrm{aug}}

\newcommand{\coeq}{\mathrm{coeq}}
\newcommand{\colim}{\mathrm{colim}}

\newcommand{\Arr}{\mathrm{Arr}}

\newcommand{\paralelas}[2]{\ar@<.2pc>[r]^-{#1} \ar@<-.2pc>[r]_-{#2}}
\newcommand{\paralelass}[2]{\ar@<.2pc>[rr]^-{#1} \ar@<-.2pc>[rr]_-{#2}}

\numberwithin{equation}{section}

\hypersetup{pdfauthor={Bruno Stonek},colorlinks=true, citecolor=blue, urlcolor=blue, linkcolor=black,breaklinks=true,hypertexnames=false}

\author{Bruno Stonek}
\address{Université Paris 13, Sorbonne Paris Cité, LAGA, CNRS, UMR 7539, 93430 Villetaneuse, France.}
\email{stonek@math.univ-paris13.fr}

\begin{document}
\title{Graded multiplications on iterated bar constructions}

\def\C{\mathcal{C}}
\newcommand{\kecalg}{k\mbox{-}\mathbf{CAlg}^\aug}
\newcommand{\kecdga}{k\mbox{-}\mathbf{cdga}^\aug}
\newcommand{\recalg}{R\mbox{-}\mathbf{CAlg}^\aug}

\begin{abstract} We define a bar construction endofunctor on the category of commutative augmented monoids $A$ of a symmetric monoidal category $\V$ endowed with a left adjoint monoidal functor $F:s\Set\to \V$. To do this, we need to carefully examine the monoidal properties of the well-known (reduced) simplicial bar construction $B_\bullet(\1,A,\1)$. We define a geometric realization $|-|$ with respect to the image under $F$ of the canonical cosimplicial simplicial set. This guarantees good monoidal properties of $|-|$: it is monoidal, and given a left adjoint monoidal functor $G:\V\to \W$, there is a monoidal transformation $|G-|\Rightarrow G|-|$. We can then consider $BA=|B_\bullet A|$ and the iterations $B^nA$. We establish the existence of a graded multiplication on these objects, provided the category $\V$ is cartesian and $A$ is a ring object. The examples studied include simplicial sets and modules, topological spaces, chain complexes and spectra.
\end{abstract}

\maketitle

The \emph{bar construction} is an algebraic machine showing up in different settings: the flavor that we consider allows for the construction, for example, of classifying spaces $BG$ of topological commutative monoids $G$, of reduced Hochschild homology $HH^k(A,k)$ where $A$ is a commutative augmented algebra over a commutative ring $k$, and of reduced topological Hochschild homology $THH^R(A,R)$ where $A$ is a commutative augmented algebra over a commutative ring spectrum $R$.

A very general incarnation of its \emph{two-sided} version appears in \cite[XII.1.1]{ekmm}: given a category $\C$, a monad $T$, a $T$-algebra $C$ and a right $T$-functor $F$, we can build a simplicial object $B_\bullet(F,T,C)$ in $\C$. A particular case is given as follows: if $\V$ is a symmetric monoidal category, $A$ is a monoid in $\V$, $N$ is a right $A$-module and $M$ is a left $A$-module, then $B_\bullet(M,A,N)\in s\V$, where $s\V$ denotes the category of simplicial objects of $\V$.

We shall focus in the following particular case, which we call the \emph{simplicial (reduced) bar construction}: it is the object of study of Section 1. Denote by $\1$ the unit of $\V$. If $A$ is an \emph{augmented} commutative monoid in $\V$, meaning that it comes with a commutative monoid morphism $A\to \1$, then $\1$ becomes a left and right $A$-module, and we can consider $B_\bullet(\1,A,\1)$. This is the object that we will denote $B_\bullet A$.

This construction is a strong symmetric monoidal functor from augmented commutative monoids in $\V$ to simplicial augmented commutative monoids in $\V$, as is well-known. A fact which to our knowledge has not appeared in the literature is the 2-categorical aspect of it of \textbf{Corollary \ref{natural}}: we explain how a strong symmetric monoidal functor $F:\V\to \W$ begets a natural isomorphism $B_\bullet F \Rightarrow FB_\bullet$.\\ 

Our primary interest lies, however, in the \emph{iteration} of the bar construction. For this we need a notion of geometric realization, from simplicial objects in $\V$ down to objects in $\V$, which we study in Section 2. This is not hard to understand abstractly. Given a fixed cosimplicial object $D$ in $\V$, the geometric realization $|X|$ of $X\in s\V$ could be defined as the tensor product of functors $X\otimes_\DDelta D$.

However, since the simplicial bar construction is a simplicial augmented commutative monoid, we need the geometric realization to take such an object to an augmented commutative monoid, in order to be able to iterate the composition $B=|B_\bullet|$. The definition just given has no reason to satisfy that. For this reason, and because all the examples that we want to consider follow this pattern, we consider only cosimplicial objects which are induced from the canonical cosimplicial simplicial set, namely 
the Yoneda embedding (see Section \ref{rempres} for what we mean by ``induced''). The geometric realization functor defined by such a cosimplicial object does indeed preserve commutative monoids: this is \textbf{Theorem \ref{realstrong}}. 

But the advantages of having chosen such a cosimplicial object do not end there. Again, we obtain a 2-categorical result: if $G:\V\to \W$ is a strong symmetric monoidal functor which is a left adjoint, then \textbf{Theorem \ref{laposta}} provides us with a \emph{monoidal} natural isomorphism $|G-|\Rightarrow G|-|$ of functors $s\V\to \W$. The fact that this natural isomorphism is monoidal is doubly important. First, it is an interesting fact per se. As a particular instance of this theorem, we obtain for example that the isomorphism of spectra \[\sip|X_\bullet|\cong |\sip X_\bullet|\]where $X$ is a simplicial topological space, is actually monoidal (Example \ref{spectro}). Second, the monoidality of this natural isomorphism implies that it begets a natural isomorphism $|GA| \cong G|A|$ of \emph{augmented commutative monoids}, if $A$ is a simplicial augmented commutative monoid, and this is crucial to our desires.\\ 

Having set this up, given an augmented commutative monoid $A$ we can define the \emph{bar construction} $BA$ as the geometric realization of the simplicial bar construction of $A$, and in Section 3 we set out to iterate it. We get a graded augmented commutative monoid $(B^nA)_{n\geq 0}$, i.e. a sequence of augmented commutative monoids. Our quest is now to find a graded multiplication \[B^nA\otimes B^m A\to B^{n+m}A.\]

To achieve this, we will need the monoid $A$ itself to have an additional multiplication. 
We need to specialize the context: instead of considering general symmetric monoidal categories, we will need to consider \emph{cartesian} ones. Such a category allows for \emph{ring objects}. Note that a ring object is in particular an augmented commutative monoid with a trivial augmentation.  
It should also be noted that from any symmetric monoidal category we can define an interesting cartesian category, the one formed by its cocommutative comonoids.

If $A$ is a ring object in a cartesian category, then the iterated bar constructions $B^*A=(B^nA)_{n\geq 0}$ is a graded ring object. This is the content of \textbf{Theorem \ref{elteo}}. Moreover, if $G$ is a left adjoint, cartesian functor between cartesian categories, then we get a natural isomorphism of graded ring objects $B^*GA\cong GB^*A$.\\ 

In the fourth and final section we unveil the examples. When we take the category of topological spaces as our cartesian category (Example \ref{topor}), then given a topological abelian group $A$, $BA$ is the bar construction introduced by Milgram \cite{milgram}. When $A$ is discrete, $B^nA$ serves as a model for an Eilenberg-Mac Lane space $K(A,n)$. Starting from a ring $R$, the graded multiplication $K(R,n)\times K(R,m)\to K(R,n+m)$ is the one found by Ravenel and Wilson in \cite{rw80}, which gives the cup product in singular cohomology with coefficients in $R$. We can also run this machine in the category of simplicial sets (Example \ref{ej1}). If $A$ is a simplicial abelian group, then $B^nA$ gives a simplicial model for a $K(A,n)$; this is well-known. But we also get that starting from a simplicial ring $R$, $B^*R$ is a graded simplicial ring, which under the geometric realization functor corresponds to the topological construction of Ravenel and Wilson. This was to be expected.

We can also work in some symmetric monoidal category of spectra (Example \ref{bravenew}). Fix $R$ a commutative ring spectrum. Then if $A$ is an augmented commutative $R$-algebra, $BA$ is the topological Hochschild homology $THH^R(A,R)$ over $R$ relative to $R$. The iterations $B^nA$ model the \emph{higher} topological Hochschild homology $THH^{R,[n]}(A,R)=S^n \odot A$, where $\odot$ denotes the tensoring over pointed topological spaces. Now, if $A$ is a ring object in cocommutative $R$-coalgebras, then we get a graded multiplication in higher $THH$ (\ref{elprod}). Moreover, if $A$ is of the form $R[S]$ where $S$ is a ring, then $THH^{R,[*]}(R[S],R) \cong R[K(S,*)]$, a natural isomorphism of graded ring objects of cocommutative $R$-coalgebras (\ref{loultimo}). Here $R[-]$ denotes $R\wedge_\S \Sigma^\infty_+(-)$.

Analogously, we can work in the symmetric monoidal category of simplicial $k$-modules, for $k$ a commutative ring (Example \ref{simpex}). Here the iterated bar constructions will rather yield Pirashvili's \cite{pirashvili} higher reduced Hochschild homology. On the other hand, Theorem \ref{elteo} is not applicable in the context of the symmetric monoidal category of non-negatively graded chain complexes of $k$-modules, the problem being that the normalized Moore functor from simplicial modules into it is not \emph{strong} monoidal (Example \ref{degeme}). Note that if $A$ is an augmented commutative differential graded algebra, then $BA$ is one as well: this construction of $BA$ coincides with the classical one introduced by Eilenberg and Mac Lane \cite[Theorem 11.1]{onthe}.\\

We will make some use of the theory of monoidal categories. We have chosen to defer this material to the appendix.

We let $\Top$ denote the cartesian closed category of compactly generated weakly Hausdorff spaces, and 
we let $k$ be a commutative ring.

For us, a \emph{spectrum} will be understood to be an $S$-module in the sense of \cite{ekmm}. A \emph{(commutative) ring spectrum} will mean a (commutative) $S$-algebra.\\

\textbf{Acknowledgments}. The author would like to thank Ignacio López Franco for the useful discussions on category theory. The contents of this article will be a part of the author's PhD dissertation at Université Paris 13. The author acknowledges support from the project ANR-16-CE40-0003 ChroK.

\section{Simplicial bar construction} \label{sect1}

Let $\V$ be a symmetric monoidal category with unit object $\1$. The symmetric monoidal category $\Mon(\V)^\aug$ has as objects the augmented monoids of $\V$, i.e. monoids $M$ with a monoid map $M\to \1$ (see \ref{augmented}). We refer to \ref{simplicial} for some background on simplicial objects.

Recall that $s\V$, the category of simplicial objects of $\V$, is symmetric monoidal with pointwise tensor product and unit $c\1$, the constant simplicial object at $\1$.

\bdefn The \emph{simplicial (reduced) bar construction} in $\Mon(\V)^\aug$ is the functor
 \[B_\bullet:\Mon(\V)^\aug\to s\V\] defined as follows. If $A\in \Mon(\V)^\aug$ with multiplication $\mu:A\otimes A\to A$, unit $\eta:\1\to A$ and augmentation $\epsilon:A\to \1$, then $B_n(A)=A^{\otimes n}$, where $A^{\otimes 0}$ means $\1$. The faces $d_i:A^{\otimes n}\to A^{\otimes n-1}$, $i=0,\dots,n$ are defined as
\[d_0=\epsilon \otimes \id^{\otimes n-1},\]\[d_i= \id^{\otimes i-1} \otimes \mu  \otimes \id^{\otimes n-i-1} \hspace{.5cm} \text{if } i=1,\dots,n-1,\]
\[d_n=\id^{\otimes n-1} \otimes \epsilon,\] 
and the degeneracies $s_i:A^{\otimes n}\to A^{\otimes n+1}$ are \[s_i=\id^{\otimes i} \otimes \eta \otimes \id^{\otimes n-i} \hspace{.5cm} \text{ for all }i=0,\dots,n.\]
\edefn

This bar construction is a particular case of the two-sided bar construction \cite[Chapter 10]{gils}: it is $B_\bullet(\1,-,\1)$ where if $A\in \Mon(\V)^\aug$ then $\1$ is viewed as a left and right $A$-module via the augmentation $\epsilon:A\to \1$.\\ 

The proof of the following proposition is elementary.

\bprop \label{strongsm} The functor $B_\bullet:\Mon(\V)^\aug\to s\V$ is strong symmetric monoidal. 
\eprop

\bobs \label{bisimpli} The symmetric monoidal category $s\V$ itself admits a monoidal bar construction \[B_\bullet: s\Mon(\V)^\aug \to s^2\V\]
which is levelwise $B_\bullet:\Mon(\V)^\aug \to s\V$. 
\eobs

\subsection{Behavior over monoidal functors}We will now provide a comparison theorem for the bar constructions in two symmetric monoidal categories $\V,\W$ related by a functor $F:\V\to \W$. We first prove a lemma that isolates the part relating purely to the monoidal structure, and then deal with the simplicial structure.

\begin{lema} \label{pseudo} Let $F:\V\to \W$ be a lax symmetric monoidal functor between symmetric monoidal categories, with structure morphism $\nabla_{A,B}:FA\otimes FB \to F(A\otimes B)$. There are lax symmetric monoidal functors and a monoidal transformation
\[\begin{tikzcd}[column sep=small, row sep=tiny]
& \W^n \ar[rd, "\otimes_n"] \ar[dd, phantom, "\ \ \Downarrow \!\nabla_{\!n}"]  \\ 
\V^n \ar[ru, "F^n"] \ar[rd, "\otimes_n"'] && \W  \\ & \V \ar[ru, "F"'] 
\end{tikzcd}\]
for every $n\geq 1$. The functors are strong if $F$ is strong.

In particular, there are monoidal transformations
\[ \begin{tikzcd}[column sep=2cm]
\V \ar[r, bend left, "F(-)^{\otimes n}", ""{below, name=a}] \ar[r, bend right, "F((-)^{\otimes n})"', ""{above, name=b}] & \W. \ar[Rightarrow,from=a,to=b,"\nabla_{\!n}"] 
\end{tikzcd}
\]
between lax symmetric monoidal functors, which are strong if $F$ is strong. 
\bprf  
The first part of this lemma can be proven by elementary means, using Mac Lane's coherence theorems \cite[VII.2, XI.2]{maclane}. We give an alternative, more sophisticated proof, essentially due to Mike Shulman \cite{shulmanmo}.

The statement that we can define $\otimes_n$ and that it is strong symmetric monoidal, and the statement that we can define $\nabla_n$ and that it is symmetric monoidal are two instances of the same result. Lack \cite[3.6]{lack} proved a coherence theorem for pseudomonoids $X$ in monoidal bicategories $\C$. By means of it, one obtains a 1-cell $\mu_n:X^n\to X$. Moreover, when $\C$ and $X$ are symmetric, then the multiplication $\mu:X\otimes X \to X$ is strong symmetric monoidal (i.e. a strong symmetric morphism of pseudomonoids), therefore $\mu_n$ is strong symmetric, too. 
This sole result encompasses both of the above results, by taking $\C$ to be the 2-category with products $\Cat$ for the first case, and by taking $\C$ to be $\mathbf{Oplax}(\mathbf{2},\Cat)$, the 2-category with products consisting of oplax functors from the interval category to $\Cat$ for the second case.

The particular case is obtained by prewhiskering with the iterated diagonal functor $\V\to \V^n$, which is strong symmetric monoidal. Indeed, whiskering preserves monoidality of natural transformations \cite[3.21, 3.24]{aguiar}.
\eprf
\end{lema}

\bprop \label{coherence} Let $F:\V\to \W$ be a normal lax symmetric monoidal functor between symmetric monoidal categories. 
Then there is a monoidal transformation
\[\xymatrix{\Mon(\V)^\aug \ar[r]^-{B_\bullet} \ar[d]_-F & s\V \ar[d]^-F \\ \Mon(\W)^\aug \ar[r]_-{B_\bullet}  & s\W \xtwocell[-1,-1]{}\omit{<1>} }\]
which is a natural isomorphism when $F$ is strong.
\bprf First of all, the fact that $F$ is \emph{normal} lax symmetric implies the existence of such an induced normal lax symmetric $F$ on \emph{augmented} monoids (\ref{moncaug}) and on simplicial objects (\ref{simplicial}).

Gathering the monoidal transformations of Lemma \ref{pseudo} together for all $n$, we obtain a monoidal transformation
\[ \begin{tikzcd}[column sep=2cm]
\V \ar[r, bend left, "\bigsqcup\limits_{n\geq 0} F(-)^{\otimes n}" pos=0.55, ""{below, name=a}] \ar[r, bend right, "\bigsqcup\limits_{n\geq 0}F((-)^{\otimes n})"' pos=0.55, ""{above, name=b}] & \Fun(\N,\W) \ar[Rightarrow,from=a,to=b,"\nabla"] 
\end{tikzcd} \]
where $\N$ is the discrete category on the natural numbers.

All there is left is to prove is that, when we start from $\Mon(\V)^\aug$, the components of $\nabla$ are really morphisms in $s\W=\Fun(\DDelta^\op,\W)$, i.e. that they are compatible with the faces and degeneracies of the simplicial bar construction.

Let $A\in \Mon(\V)^\aug$ with multiplication $\mu:A\otimes A\to A$. The face map $d_1:B_3(A)\to B_2(A)$ is \[\mu \otimes \id:A\otimes A\otimes A\to A\otimes A,\]
and the face map $d_1:B_3(FA)\to B_2(FA)$ is the composition of the two vertical maps on the left of the following diagram, whose commutativity we need to check.
\[
\begin{tikzcd}
FA\otimes FA \otimes FA \ar[r, "\nabla_3"] \ar[d, "\nabla\otimes\id"'] & F(A\otimes A \otimes A) \ar[dd, "F(\mu\otimes \id)"] \\ F(A\otimes A)\otimes FA \ar[d, "F\mu\otimes \id"'] \\ FA\otimes FA  \ar[r, "\nabla"'] & F(A\otimes A)
\end{tikzcd}
\]
But the proof of Lemma \ref{pseudo} guarantees that we can take $\nabla_3=\nabla\circ (\nabla\otimes \id)$, therefore this commutativity is an application of naturality of $\nabla$. All the other differentials \linebreak $d_i:A^{\otimes n}\to A^{\otimes (n-1)}$ for $i\neq0,n$ at each level are built in the same fashion, so the proof adapts. For the extremal differentials $d_0$ and $d_n$ which use the augmentation $\epsilon:A\to \1$ of $A$, there is another diagram chase. For example, for $d_0:A\otimes A\otimes A \to A\otimes A$, the verification is the commutativity of the following outer diagram. The inner diagrams commute by naturality and unitality of $\nabla$.
\[
\begin{tikzcd}
FA\otimes FA \otimes FA \ar[r, "\nabla \otimes \id"] \ar[d, "F\epsilon \otimes \id \otimes \id"'] & F(A\otimes A) \otimes FA \ar[d,"F(\epsilon \otimes \id)\otimes \id"] \ar[r,"\nabla"] & F(A\otimes A\otimes A) \ar[d, "F(\epsilon \otimes \id\otimes \id)"] \\
F\1 \otimes FA \otimes FA \ar[r, "\nabla\otimes \id"] & F(\1\otimes A)\otimes FA \ar[ldd,"\cong"] \ar[r,"\nabla"] & F(\1\otimes A \otimes A) \ar[dd,"\cong"] \\ \1\otimes FA \otimes FA \ar[u,"\cong"] \\
FA\otimes FA \ar[rr,"\nabla"'] \ar[u,"\cong"] && F(A\otimes A) 
\end{tikzcd}
\]
Compatibility with the degeneracies follows from a similar diagram chase.
\eprf
\eprop

\subsection{Commutative augmented monoids}We denote by $\CMon(\V)$ the symmetric monoidal category of commutative monoids in $\V$ (see \ref{general}), and by $\Ab(\C)$ the cartesian category of abelian group objects in a cartesian category $\C$ (see \ref{cartesiansect}). Combining Propositions \ref{strongsm} and \ref{coherence} we obtain the following result.

\bcor \label{natural} There is an induced strong symmetric monoidal functor \[B_\bullet: \CMon(\V)^\aug\to s(\CMon(\V)^\aug).\]
Let $F:\V\to \W$ be a normal lax symmetric monoidal functor between symmetric monoidal categories. There is a  
natural transformation 
displayed in the following diagram, which is a natural  isomorphism when $F$ is strong.
\[\xymatrix{\CMon(\V)^\aug \ar[r]^-{B_\bullet} \ar[d]_-F & s(\CMon(\V)^\aug) \ar[d]^-F \\ \CMon(\W)^\aug \ar[r]_-{B_\bullet} & s(\CMon(\W)^\aug) \xtwocell[-1,-1]{}\omit{<0>} }\]
If $\C$ is a cartesian category, there is an induced cartesian functor \[B_\bullet:\Ab(\C)\to s\Ab(\C),\] and if $F:\C\to \D$ is a cartesian functor between cartesian categories, then there is a natural isomorphism
\[\xymatrix{\Ab(\C) \ar[r]^-{B_\bullet} \ar[d]_-F & s\Ab(\C) \ar[d]^-F \\ \Ab(\D) \ar[r]_-{B_\bullet} & s\Ab(\D). \xtwocell[-1,-1]{}\omit{<0>} }\]
\bprf  The first claim follows from Proposition \ref{strongsm} by passing to commutative augmented monoids, and applying the Eckmann-Hilton argument (\ref{general}) and the arguments in the appendix which allow us to permute the order of taking simplicial objects and augmented objects. 

The second claim follows from applying $\CMon(-)^\aug$ to the diagram in Proposition \ref{coherence} and using the result from \ref{moncaug} which states that a monoidal transformation between normal lax symmetric monoidal functors indeed induces a natural transformation in the respective categories of augmented commutative monoids.  

If $\C$ is cartesian, then $\Ab(\C)$ and $s\C$ are also cartesian, and Proposition \ref{strongsm} says that $B_\bullet:\Mon(\C)^\aug\to s\C$ is cartesian. The proof is now as above, but passing to abelian group objects instead of augmented commutative objects. Indeed, there is an equivalence of categories $\Ab(\C)^\aug\cong \Ab(\C)$ (\ref{moncaug}).
\eprf
\ecor

\section{Geometric realization} 

We now need a notion of ``intrinsic'' geometric realization of a simplicial object in a category $\V$ to an object in $\V$. By the quoted term we mean a functor $s\V\to \V$, in contrast with what happens with the standard geometric realization of a simplicial set into a topological space (see Example \ref{topor}), which we call ``extrinsic'' and which we analyze in Remark \ref{enrich}.\ref{ajoba}. Recall that the geometric realization of a simplicial space can be described as a tensor product of functors \cite[IX.6]{maclane}: if $X_\bullet$ is a simplicial space, then
\[|X_\bullet|=X_\bullet\otimes_\DDelta \Delta^\bullet_{\textup{top}} \in \Top\]
where $\Delta^\bullet_{\textup{top}}:\DDelta\to \Top$ is the standard cosimplicial space that takes $[n]$ to the topological $n$-simplex. This defines a functor $|-|:s\Top\to \Top$. 

More generally, we have the following definition.

\bdefn \label{def31}  Let $\V$ be a cocomplete symmetric monoidal category with a given cosimplicial object $D^\bullet:\DDelta \to \V$. We define the \emph{geometric realization (with respect to $D^\bullet$)} as the functor
\begin{equation}\label{ten1} |-|_{D^\bullet}\coloneqq -\otimes_\DDelta D^\bullet: s\V\to \V.\end{equation}
Explicitly, if $X_\bullet \in s\V$, then $|X_\bullet|_{D^\bullet}$ can be expressed as a coend, or even more explicitly as a coequalizer, as follows: 
\[|X_\bullet|_{D^\bullet}= \int^{n} X_n \otimes D^n = \coeq \left( \bigsqcup_{[n]\stackrel{f}{\to} [m] \in \Arr \DDelta} X_m \otimes D^n \rightrightarrows \bigsqcup_{[p]\in \DDelta} X_p \otimes D^p \right)\]
where the parallel arrows are induced by the maps $X_m\otimes D^n \to X_n \otimes D^n$ and $X_m \otimes D^n \to X_m \otimes D^m$ given by the action of $f$. 
\edefn

\begin{lema} \label{adjunto}
If $\V$ is moreover \emph{closed}, then the functor $|-|_{D^\bullet}:s\V\to \V$ is a left adjoint. The right adjoint is given by the functor
\[\underline{\V}(D^\bullet,-):\V\to s\V, \hspace{.5cm} V\mapsto ([n]\mapsto \underline{\V}(D^n,V))\]
where $\underline{\V}(-,-)\in \V$ denotes the internal hom object of $\V$.
\bprf The proof is an elementary manipulation of adjunctions and coends, which uses that the set $\Nat(F,G)$ of natural transformations from $F$ to $G$ is equal to the end $\int_X \V(FX,GX)$ for functors $F$ and $G$ out of some small category.
\eprf
\end{lema}

\bcor \label{const} The geometric realization of a simplicial object $cX\in s\V$ constant at $X\in \V$ is isomorphic to $X\otimes D^0$.
\bprf
Indeed,
\begin{align*}
\V(|cX|_{D^\bullet},V) &\cong \Nat(cX,\underline{\V}(D^\bullet,V)) \cong \V(X,\lim \underline{\V}(D^\bullet,V) ) 
\cong \lim \V(X,\underline{\V}(D^\bullet,V)) \cong \\
&\cong \lim \V(X\otimes D^\bullet,V) \cong \V(\colim (X\otimes D^\bullet),V) 
\cong \V(X\otimes D^0,V)
\end{align*}
naturally in $V\in \V$. Yoneda's lemma implies $|cX|_{D^\bullet}\cong X\otimes D^0$. Here we have used the fact that the colimit of a cosimplicial object is its zeroth component, since $\DDelta$ has $[0]$ as its final object. 
\eprf
\ecor

\bobs \label{enrich} It is interesting to note that geometric realization can be extended to an enriched context, in two different ways. Let $\C$ be a cocomplete category enriched, tensored and cotensored over a closed symmetric monoidal category $\V$.
\be
\item  If we are given a cosimplicial object $D^\bullet:\DDelta \to \V$, then we can define $|-|_{D^\bullet}:s\C\to \C$ as $|X_\bullet|_{D^\bullet}=\int^n X_n \odot D^n$ where $\odot$ denotes the tensoring of $\C$ over $\V$. The proof of Lemma \ref{adjunto} adapts to prove that $|-|_{D^\bullet}:s\C\to \C$ is a left adjoint. If we take $\C=\V$, we recover the construction above.
\item \label{ajoba} If instead we are given a cosimplicial object $D^\bullet:\DDelta \to \C$, then we can define $|-|_{D^\bullet}^e:s\V\to \C$ as $|X_\bullet|_{D^\bullet}^e=\int^n D^n \odot X_n$. Again, $|-|^e_{D^\bullet}$ is a left adjoint, and taking $\C=\V$ recovers the construction above. Since this functor takes simplicial objects in one category and begets objects in a different category, we call it \emph{extrinsic geometric realization}. 
\ee
\eobs

\subsection{Cosimplicial objects induced by the Yoneda embedding} \label{rempres}

We will concentrate on closed symmetric monoidal categories with a cosimplicial object induced by the Yoneda embedding on simplicial sets. More precisely, 

\bq From here on we let $\V$ be a cocomplete, closed symmetric monoidal category together with a 
lax symmetric monoidal functor $F:s\Set\to \V$ which is a left adjoint. Let $\Delta^\bullet: \DDelta \to s\Set$ be the Yoneda embedding: we consider the cosimplicial object on $\V$ given by $F\Delta^\bullet:\DDelta \to \V$. We denote by $|-|$ the geometric realization functor $|-|_{F\Delta^\bullet}:s\V\to \V$.
\eq

An additional hypothesis of normality for $F$ will be specified when necessary; notably, whenever we need $F$ to induce a functor on \emph{augmented} objects. 

The particularities of the Yoneda embedding will allow us to do things we could not do with the geometric realization with respect to an abstract cosimplicial object (see e.g. formula (\ref{kan})).

\bobs \label{const2} Let $X\in \V$. Since $\Delta^0$ is the unit of the cartesian category $s\Set$, the unit of $F$ has the form $\nabla_0:\1 \to F\Delta^0$. Therefore we get an arrow
\begin{equation}\label{const-arrow} \xymatrix@C+1pc{X \ar[r]^-{\rho^{-1}}_-\cong & X\otimes \1 \ar[r]^-{\id \otimes \nabla_0} & X \otimes F\Delta^0 \cong |cX| } \end{equation}
which is an isomorphism if $F$ is normal. The last isomorphism is provided by Corollary \ref{const}.
\eobs

We recall the \emph{density theorem} \cite[1.4.6]{riehl}: if $H:I^\op\to \C$ is a functor from the opposite of a small category $I$ into a cocomplete category $\C$, then 
\begin{equation}\label{density}Hj \cong \int^i I(j,i)\cdot Hi.\end{equation}
Here and henceforth $\cdot$ denotes the tensoring of $\C$ over $\Set$. Explicitly, if $A\in \Set$ and $C\in \C$, then $A\cdot C$ is the coproduct $\bigsqcup_{i\in A} C$.

Note that since $F:s\Set \to \V$ is a left adjoint it preserves coproducts, and therefore $F(A\cdot X)\cong A\cdot FX$.

\bteo  \label{realstrong} The geometric realization functor $|-|:s\V\to \V$ has a lax symmetric monoidal structure, which is strong (resp. normal) if $F$ is strong (resp. normal).
\bprf Let us express the simplicial set $\Delta^n\times \Delta^m$ as a coend, using the density theorem.
\begin{align} \label{kan}
\nonumber (\Delta^n \times \Delta^m)(j)&\cong \int^i \DDelta(j,i) \times (\Delta^n \times \Delta^m)(i)= \int^i \Delta^i(j) \times (\DDelta \times \DDelta)((i,i),(n,m)) \\
&= \left(\int^i (\DDelta \times \DDelta)((i,i),(n,m)) \cdot \Delta^i\right)(j)
\end{align}
Let $X_\bullet,Y_\bullet\in s\V$. We will now repeatedly use Fubini's theorem for coends \cite[IX.8]{maclane}, the fact that $-\otimes -$ commutes with colimits (hence with coends) separately in each variable, and the fact that $F$ commutes with coends. Finally, we use the density theorem again.
\begin{align*}
|X_\bullet|\otimes |Y_\bullet| &= \int^n X_n\otimes F\Delta^n \otimes \int^m Y_m \otimes F\Delta^m \\
&\cong \int^{n,m} X_n \otimes Y_m \otimes F\Delta^n\otimes F\Delta^m \\
&\to\int^{n,m} X_n\otimes Y_m \otimes F(\Delta^n\times \Delta^m) \\
&\cong \int^{n,m} X_n\otimes Y_m \otimes F\left( \int^i (\DDelta \times \DDelta)((i,i),(n,m)) \cdot \Delta^i \right) \\
&\cong \int^{n,m} X_n \otimes Y_m \otimes \int^i (\DDelta \times \DDelta)((i,i),(n,m)) \cdot F\Delta^i \\
&\cong \int^i \left( \int^{n,m} (\DDelta \times \DDelta)((i,i),(n,m)) \cdot (X_n\otimes Y_m) \right) \otimes F\Delta^i \\
&\cong \int^i X_i\otimes Y_i \otimes F\Delta^i = |X_\bullet\otimes Y_\bullet|
\end{align*}
The unit of $|-|$ is furnished by (\ref{const-arrow}) applied to $X=\1$.

It is a long, if tedious verification that these morphisms endow $|-|$ with the structure of a lax symmetric monoidal functor. 
\eprf
\eteo

\subsection{Behavior under monoidal functors}

\bq From here on we let $\W$ be a cocomplete, closed symmetric monoidal category together with a 
lax symmetric monoidal functor $G:\V\to \W$ which is a left adjoint. We endow $\W$ with the cosimplicial object $GF\Delta^\bullet:\DDelta \to \W$, and we denote by $|-|$ the geometric realization functor $|-|_{GF\Delta^\bullet}:s\W\to \W$.\eq

Of course, Theorem \ref{realstrong} applies to $\W$ as well: $|-|:s\W\to \W$ has a 
lax symmetric monoidal structure which is strong (resp. normal) if $F$ and $G$ are strong (resp. normal). Note that we will be using $|-|$ to denote both $|-|_{F\Delta^\bullet}$ and $|-|_{GF\Delta^\bullet}$: which one we mean should be clear from context.

\bteo  \label{laposta} There is a monoidal transformation between  
lax symmetric monoidal functors
\begin{equation}\label{presomon} \xymatrix@C+1pc{s\V \rtwocell<5>^{|G-|}_{G|-|} {\tau}  & \W}\end{equation}
which is an isomorphism if $G$ is strong.
\eteo
In particular, this holds when $\V=s\Set$ and $F=\id_{s\Set}$.

\bprf Since $G$ is a left adjoint, it preserves coends, and thus we get
\[ \int^n GX_n \otimes GF\Delta^n  \to \int^n G(X_n\otimes F\Delta^n) \cong G\int^n X_n \otimes F\Delta^n \]
for $X_\bullet\in s\V$, defining the desired natural transformation. We need to check it is monoidal, i.e. that the following diagram commutes, for $X_\bullet, Y_\bullet\in s\V$.

\[\xymatrix{|GX_\bullet| \otimes |GY_\bullet|\ar[r] \ar[d]_-{\tau\otimes \tau} & |GX_\bullet\otimes GY_\bullet| \ar[r] & |G(X_\bullet \otimes Y_\bullet)| \ar[d]^-\tau \\ G|X_\bullet| \otimes G|Y_\bullet| \ar[r] & G(|X_\bullet|\otimes |Y_\bullet|) \ar[r] & G|X_\bullet \otimes Y_\bullet| }\]
The horizontal arrows intertwine the monoidal structure of $G$ and of its induced functor $G:s\V\to s\W$ with the monoidal structure of the geometric realizations $|-|$ in $s\V$ and in $s\W$ obtained in Theorem \ref{realstrong}. These structures are defined via a fairly long string of isomorphisms, whence the difficulty of reproducing the necessary diagram proof. One can expand the diagram into one big rectangle filled with coends, and juggle around with naturality and monoidality properties of $F$ and $G$. The gist of the proof is that the geometric realizations in $\V$ and in $\W$ are not independent: they are related via $G$, which is 
lax symmetric monoidal. 
\eprf

\subsection{Commutative augmented monoids} Suppose $F$ and $G$ are normal. 
By Theorem \ref{realstrong}, the geometric realization functor $|-|:s\V\to \V$ is normal lax symmetric monoidal, therefore it induces a functor\begin{equation}\label{gro}|-|:s\CMon(\V)^\aug \to \CMon(\V)^\aug\end{equation}
and similarly for $|-|:s\W\to \W$.

Since the natural transformation of Theorem \ref{laposta} is monoidal, we can apply (\ref{cosados}) and obtain a natural transformation
\begin{equation}\label{presomon2} \xymatrix@C+1pc{s\CMon(\V)^\aug \rtwocell<6>^{|G-|}_{G|-|} {\tau}  & \CMon(\W)^\aug}\end{equation}
which is an isomorphism if $G$ is strong.

\section{The bar construction}

Suppose $F$ is normal. 
We can now define the \emph{bar construction}
\begin{equation} \label{bolein} B:\CMon(\V)^\aug\to \CMon(\V)^\aug\end{equation}
as the composite $\xymatrix{\CMon(\V)^\aug\ar[r]^-{B_\bullet} & s(\CMon(\V))^\aug \ar[r]^-{|-|} & \CMon(\V)^\aug}$. 
Remark that the monoid structure on $BA$ for $A\in \CMon(\V)^\aug$ is induced by the simplicial map \begin{equation}\label{mas}A^{\otimes p}\otimes A^{\otimes p}\to A^{\otimes p}\end{equation} which is the monoid structure on $A^{\otimes p}$.

If $\V=\C$ is cartesian and $F$ is cartesian, we can take abelian group objects instead of augmented commutative monoids, and obtain a functor \[B:\Ab(\C)\to \Ab(\C).\] 

\bprop \label{natural2} 
Suppose $F$ and $G$ are normal. 
There is a natural transformation 
\[\xymatrix{ \CMon(\V)^\aug \ar[r]^-B \ar[d]_-G & \CMon(\V)^\aug \ar[d]^-G \\ \CMon(\W)^\aug \ar[r]_-B & \CMon(\W)^\aug \xtwocell[-1,-1]{}\omit{<0>}.  }\]
which is an isomorphism if $G$ is strong.

If $\V$, $\W$, $F$ and $G$ are cartesian, 
there is an analogous square for $\Ab$ instead of $\CMon^\aug$.
\eprop

\bprf The 
natural transformation is the pasting of the following two.
\[\xymatrix{\CMon(\V)^\aug \ar[d]_-G \ar[r]^-{B_\bullet} & s\CMon(\V)^\aug \ar[d]^-G \ar[r]^-{|-|} & \CMon(\V)^\aug \ar[d]^-G  \\
\CMon(\W)^\aug \ar[r]_-{B_\bullet} & s\CMon(\W)^\aug \ar[r]_-{|-|} \xtwocell[-1,-1]{}\omit{<0>} & \CMon(\W)^\aug \xtwocell[-1,-1]{}\omit{<0>} }\]
The left one comes from Corollary \ref{natural}, and the right one is (\ref{presomon2}). 
For abelian groups it is entirely analogous. 
\eprf

\subsection{Remark: bicommutative Hopf monoids} \label{remark:bihopf}

By the results of \ref{cartesiansect}, the category $\CoComon(\V)$ of cocommutative comonoids in $\V$ is cartesian, therefore admits a simplicial bar construction on its abelian group objects (Corollary \ref{natural}). As remarked in \ref{hopfmo}, $\Ab(\CoComon(\V))$ is equivalent to $\biHopf(\V)$, the category of bicommutative Hopf monoids in $\V$, so the simplicial bar construction in this case is actually a functor \linebreak 
$B_\bullet: \biHopf(\V) \to s(\biHopf(\V))$. This functor coincides with the  simplicial bar construction of Corollary \ref{natural} under the forgetful functors $U$ from (simplicial) bicommutative Hopf monoids to (simplicial) commutative augmented monoids.

Suppose $F:s\Set\to \V$ is strong. As seen in Theorem \ref{realstrong}, this implies that $|-|:s\V\to \V$ is strong symmetric monoidal. In the same fashion as above, we obtain an induced functor $|-|:s\biHopf(\V)\to \biHopf(\V)$ which coincides with (\ref{gro}) under the forgetful functors $U$.

Composing the functors we get a functor \begin{equation}\label{barquito} B:\biHopf(\V)\to \biHopf(\V)\end{equation} which coincides with (\ref{bolein}) under $U$.

\subsection{Iterated bar constructions and ring structure}

We can iterate the bar construction to obtain functors \[B^n:\CMon(\V)^\aug\to \CMon(\V)^\aug \hspace{.3cm} \text{ for } n\geq 0,\]
where we define $B^0$ to be the identity functor.

If $\V=\C$ is cartesian, we similarly obtain functors
\[B^n:\Ab(\C)\to \Ab(\C) \hspace{.3cm} \text{ for } n\geq 0.\]

We will now put a graded ring structure on iterated bar constructions, provided we start with a ring object. To be able to carry this out, we make the following assumptions.

\bq We let the monoidal structures in our categories $\V, \W$ be cartesian. We therefore change notation: we are given cartesian categories $\C, \D$, a cartesian functor $F:s\Set \to \C$ and a cartesian functor $G:\C\to \D$, both of which are left adjoints. \eq 

Thus we get induced geometric realizations on $s\C$ and on $s\D$ which are cartesian.

We can glue all the iterated bar constructions together into a single functor \[B^*:\Ab(\C)\to \NGrAb(\C), \hspace{.5cm} A\mapsto (B^nA)_{n\in \N}.\]

Here $\NGrAb(\C)$ stands for the category of $\N$-graded objects of $\Ab(\C)$, i.e. the functor category $\Fun(\N,\Ab(\C))$. 

\bteo \phantomsection \label{elteo} \be[a)] \item The functor $B^*$ extends to a functor \[B^*:\Ring(\C)\to \NGrRing(\C).\]
\item \label{elteob} There is a natural isomorphism
\[\xymatrix{ \Ring(\C) \ar[d]_-G \ar[r]^-{B^*} & \NGrRing(\C) \ar[d]^-{G} \\ \Ring(\D) \ar[r]_-{B^*} & \NGrRing(\D). \xtwocell[-1,-1]{}\omit{<0>} }\]
\ee
\eteo
Here $\NGrRing(\C)$ stands for the category of graded ring objects in $\C$: its objects are sequences $(A_n)_{n\geq 0} \in \NGrAb(\C)$ together with a graded multiplication $A_n\times A_m \to A_{n+m}$ which is associative, unital and distributive with respect to the abelian group operation of each $A_i$.
\bprf \be[a)] \item 
Let $S\in \Ring(\C)$. Denote by $\mu:S\times S \to S$ its multiplication. We define the graded multiplication $\smile_{n,m}:B^nS \times B^mS \to B^{n+m}S$ inductively.

For $n=m=0$, $\smile_{0,0}:S\times S\to S$ is $\mu$. Now let us define $\smile_{0,m+1}$ from $\smile_{0,m}$.

Consider, for $i=1,\dots,p$, the composition \begin{equation}\label{eq}\xymatrix@C+1pc{S\times (B^mS)^{\times p} \ar[r]^-{\id \times \pi_i} & S\times B^mS \ar[r]^-{\smile_{0,m}} & B^mS,} \end{equation}
where $\pi_i$ denotes the $i$-th projection map. By the universal property of the categorical product, these maps define a morphism
\[\varphi_m^p: S \times (B^mS)^{\times p} \to (B^mS)^{\times p}\]in $\C$ which commutes with the faces and degeneracies of $B_\bullet B^mS$. Indeed, as an example, the commutativity with the differentials $d_1,\dots,d_{p-1}$ rests on the distributivity of $\smile_{0,m}$ with respect to the addition of $B^mS$, and this is obtained inductively: for $m=0$ it is the mere distributivity of $\smile_{0,0}\ =\mu$ with respect to addition. We say more about the distributivity of the higher $\smile_{n,m}$ at the end of the proof.

We thus get a morphism
\[\varphi_m: S\times B_\bullet B^mS\to B_\bullet B^mS\]
in $s\C$, where $S$ is seen as a constant simplicial object.

As the geometric realization of a constant simplicial object gives the original object (Remark \ref{const2}) and as geometric realization is a cartesian functor, we obtain an induced map 
\begin{equation} \label{komar} \xymatrix{S\times B^{m+1}S \ar[r]^-\cong & |S\times B_\bullet B^mS| \ar[r]^-{|\varphi_m|} & B^{m+1}S}\end{equation}
which we call $\smile_{0,m+1}$.

The definition of $\smile_{n+1,m}$ from $\smile_{n,m}$ is symmetrical: replace (\ref{eq}) with \[\xymatrix@C+1pc{(B^nS)^{\times p} \times B^mS \ar[r]^-{\pi_i \times \id} & B^nS\times B^mS \ar[r]^-{\smile_{n,m}} & B^{n+m}S}\]
and repeat the process.

The unit for this graded multiplication is the unit for the multiplication of $S$. Associativity and distributivity of $\smile$ follow from associativity and distributivity of $\mu$; these are all straightforward verifications. As an example, here is the diagram expressing the distributivity of $\smile_{0,1}$ at the simplicial level,
\[\xymatrix@C+1pc{S\times S^{\times p} \times S^{\times p} \ar[r]^-{\Delta \times \id \times \id} \ar[dd]_-{\id \times +} & S \times S \times S^{\times p} \times S^{\times p} \ar[r]^-{\id \times \sigma \times \id} & S \times S^{\times p} \times S \times S^{\times p} \ar[d]^-{\varphi_0^p \times \varphi_0^p} \\ && S^{\times p} \times S^{\times p} \ar[d]^-+ \\ S \times S^{\times p} \ar[rr]_-{\varphi_0^p} && S^{\times p}} \]
where $+:S^{\times p} \times S^{\times p} \to S^{\times p}$ is the abelian group structure map of $S^{\times p}$ (\ref{mas}), $\Delta:S\to S\times S$ is the diagonal, and $\sigma$ is the symmetry. Its commutativity follows from distributivity in $S$.

\item First of all, the vertical functors induced by $G$ do exist since $G:\C\to \D$ is cartesian. Second, the commutativity at the level of each abelian group object follows by iterating Proposition \ref{natural2}.

For the compatibility of $\smile$-multiplications, first observe that when $n=m=0$ this is just the definition of the multiplication in $GS$, for $S\in \Ring(\C)$:
\[\xymatrix{GS \times GS \ar[r]^-{\mu^{GS}} \ar[d]_-{\cong} & GS \ar@{=}[d] \\ G(S\times S) \ar[r]_-{G\mu^S} & GS.}  \]
For general $n$ and $m$, this amounts to the commutativity of the following diagram in $\D$,
\[\xymatrix@C+1pc{B^nGS \times B^mGS \ar[r]^-{\smile_{n,m}^{GS}} \ar[d]_-\cong & B^{n+m}GS \ar[dd]^-\cong \\ GB^nS \times GB^mS \ar[d]_-{\cong} \\ G(B^nS \times B^mS)\ar[r]_-{G(\smile_{n,m}^S)} & G(B^{n+m}S)}\]
which holds since the definition of the $\smile$-multiplications only involves arrows which commute with $G$ and $B$. \qedhere
\ee
\eprf 

\bobs We really need that $G$ be cartesian, since we need it to preserve ring objects, but we could ask that $F$ be merely normal lax symmetric monoidal. This affects the proof only in that the isomorphism in (\ref{komar}) becomes just a morphism.
\eobs

\subsection{Cocommutative comonoids} \label{cocaca}

If our categories are symmetric monoidal but not cartesian, we cannot a priori carry out the construction of the previous section. However, note that $s\Set$ \emph{is} cartesian, therefore for any strong symmetric monoidal functor $F:s\Set\to \V$ we obtain a cartesian functor between cartesian categories
\[F:s\Set \to \CoComon(\V)\]
as remarked in (\ref{laec}). Moreover, if $G:\V\to \W$ is strong symmetric monoidal, we also obtain a cartesian functor
\[G:\CoComon(\V)\to\CoComon(\W).\]
Thus, from Theorem \ref{elteo} we obtain the following result.
\bcor \phantomsection \label{elteocom}
\be
\item Let $\V$ be a cocomplete closed symmetric monoidal category and $F:s\Set\to \V$ be a strong symmetric monoidal functor which is a left adjoint. There is an iterated bar construction functor 
\begin{equation}\label{to}B^*: \Ring(\CoComon(\V)) \to \NGrRing(\CoComon(\V)).\end{equation}
\item Let $\W$ be a cocomplete closed symmetric monoidal category and $G:\V\to \W$ be a strong symmetric monoidal functor which is a left adjoint. There is a natural isomorphism
\[\xymatrix{ \Ring(\CoComon(\V)) \ar[d]_-G \ar[r]^-{B^*} & \NGrRing(\CoComon(\V)) \ar[d]^-{G} \\ \Ring(\CoComon(\W)) \ar[r]_-{B^*} & \NGrRing(\CoComon(\W)). \xtwocell[-1,-1]{}\omit{<0>} }\]
\ee
\ecor
 
\bdefn \label{ladef} A graded ring object in cocommutative comonoids in $\V$ is called a \emph{coalgebraic} (or \emph{Hopf}) \emph{ring in $\V$}. We call a coalgebraic ring in $\kmod$ a $k$\emph{-coalgebraic ring}. 
\edefn

The notion of $k$-coalgebraic ring was introduced in \cite{rw77}: they called it ``Hopf ring''. We prefer this other term (which they also considered but didn't keep) since it is more explicit.\\

We can apply the second part of this Corollary to a strong symmetric monoidal, left adjoint functor $F:s\Set\to \V$, in which case the situation simplifies: there is a natural isomorphism
\begin{equation}\label{elteocoms} \xymatrix{ s\Ring \ar[d]_-F \ar[r]^-{B^*} & s\NGrRing \ar[d]^-{F} \\ \Ring(\CoComon(\V)) \ar[r]_-{B^*} & \NGrRing(\CoComon(\V)). \xtwocell[-1,-1]{}\omit{<0>} }\end{equation}

\bobs \label{nojodas} There is a forgetful functor \[\Ring(\CoComon(\V))\to \Ab(\CoComon(\V))=\biHopf(\V).\]Thus, the functor (\ref{to}) has each of its levels $B^n$ forget down to the respective iteration of (\ref{barquito}).

\eobs

\section{Examples}
\subsection{Simplicial sets} \label{ej1}We start with the cartesian closed category $s\Set$ itself. Endowed with the Yoneda embedding $\Delta^\bullet$ as a cosimplicial simplicial set, the induced geometric realization \[|-|_{\Delta^\bullet}:s^2\Set\to s\Set\] is naturally isomorphic to the diagonal functor, as is well-known. In symbols,
\[\int^n X_{n,\bullet} \times \Delta^n \cong \diag(X).\]
A proof for this formula can be obtained by means of the density theorem (\ref{density}), the coend formula for $\Delta^n\times \Delta^m$ (\ref{kan}) and Fubini's theorem for coends.

We obtain a functor \begin{equation}\label{bsset}B:s\Ab\to s\Ab.\end{equation}

It is a classical result that this functor is weakly homotopy equivalent to the $\bar{W}$-construction of Eilenberg and Mac Lane \cite{onthe}. Some discussion and references for this can be found in \cite{stevenson} after Lemma 15. One way of obtaining this result is as follows. Duskin identified the functor $\bar{W}:s\Ab\to s\Ab$ with the functor $TB_\bullet:s\Ab\to s\Ab$, where $T:s^2\Ab\to s\Ab$ is a functor going by several names, two of which are ``Artin-Mazur diagonal'' and ``totalization''. Then one needs to provide a natural weak homotopy equivalence $T\Rightarrow \diag$. This result has a complicated history: we refer to the aforementioned discussion by Stevenson, and to the more recent \cite{zisman}.  

If $G$ is an abelian group and we view it as a constant simplicial abelian group, then Eilenberg and Mac Lane proved that $\bar{W}^nG$ gives a simplicial abelian group model for a $K(G,n)$; see \cite[A.21]{dtp}. 

Theorem \ref{elteo} gives us a functor $B^*:\Ring(s\Set)\to \NGrRing(s\Set)$, i.e.
\begin{equation}\label{bring} B^*:s\Ring \to s\Gr\Ring. \end{equation} When $S$ is a constant simplicial ring, $B^nS$ is a simplicial model for $K(S,n)$, and the graded multiplication is a simplicial model for the cup product in Eilenberg-Mac Lane simplicial sets, as we will see in the next section.

\subsection{Topological spaces} \label{topor}
Recall that $\Top$ denotes the cartesian closed category of compactly generated weakly Hausdorff spaces. Let \[F=|-|^e:s\Set\to \Top\] be the classical ``extrinsic'' geometric realization functor. It is cartesian: this is a well-known result of Milnor on the geometric realization of a product of simplicial sets. It is also a left adjoint: its right adjoint is the singular functor.

Note that this extrinsic geometric realization functor is defined following the pattern of Remark \ref{enrich}.\ref{ajoba}: the base (cartesian) monoidal category is $\Set$, and we consider $\Top$ as a $\Set$-category. 

The cosimplicial object $|\Delta^\bullet|^e$ is the standard cosimplicial space, i.e. $|\Delta^n|^e$ is the topological $n$-simplex. Therefore the geometric realization $|-|_{|\Delta^\bullet|^e}:s\Top \to \Top$ is the standard geometric realization of a simplicial space, as considered e.g. in \cite{gils}. It should be noted that our Theorem \ref{realstrong} gives a categorical proof that it preserves products (compare with \cite[11.5]{gils}): the topology is contained entirely in Milnor's theorem that $|-|^e$ is cartesian.

The resulting functor \[B:\Ab(\Top)\to \Ab(\Top)\] is Milgram's \cite{milgram} version of the bar construction of a topological abelian group, as observed by Mac Lane \cite{maclanemilgram}. The space $BG$ is an especially nice model for the classifying space of $G$, because it carries a strict topological abelian group structure. Thus if $G\in \Ab(\Top)$ is discrete, then $BG$ is a model for the Eilenberg-Mac Lane space $K(G,1)$, and $B^nG$ is a model for a $K(G,n)$.

Theorem \ref{elteo} applied to $F=\id_{s\Set}$ and $G=|-|^e:s\Set \to \Top$ gives a natural isomorphism
\[\xymatrix{ s\Ring \ar[r]^-{B^*} \ar[d]_-{|-|^e} & \NGrRing(s\Set) \ar[d]^-{|-|^e} \\ \Ring(\Top) \ar[r]_-{B^*} & \NGrRing(\Top). \xtwocell[-1,-1]{}\omit{<0>}  }\]
In other words, we have natural isomorphisms \begin{equation}\label{hojo}B^n(|S|^e)\cong |B^nS|^e\end{equation} compatible with the graded multiplications existing on each side, for $S\in s\Ring$.

When $S\in \Ring(\Top)$ is a discrete topological ring, $B^*S=(K(S,n))_{n\geq 0}$ where $K(S,n)$ denotes an $n$-th Eilenberg-Mac Lane space of $S$, and the graded multiplication \begin{equation}\label{cupo}\smile:K(S,n)\times K(S,m)\to K(S,n+m)\end{equation}  represents the cup product in ordinary cohomology with coefficients in $S$ \cite[1.7]{rw80}. 
Thus, viewing $S$ as a constant simplicial ring, the graded multiplication in simplicial sets
\begin{equation}\label{barsim} B^nS \times B^m S \to B^{n+m}S\end{equation}
coincides with the cup product map (\ref{cupo}) after geometric realization under the isomorphism (\ref{hojo}), i.e. we have gotten a simplicial construction of the cup product map in Eilenberg-Mac Lane simplicial sets.\\

Let us pass to homology. Let $E$ be a commutative ring spectrum. We denote by $E_*=\pi_*(E)$ its graded commutative ring of coefficients. Let $E_*(-):\Top \to E_*\Mod$ denote its associated unreduced homology theory on spaces taking values in $E_*$-graded modules.

The category $E_*\Mod$ is symmetric monoidal with the tensor product $\otimes_{E_*}$. The functor $E_*(-)$ has a lax symmetric monoidal structure given by the homological cross product
\begin{equation}\label{kun}E_*(X)\otimes_{E_*} E_*(Y)\to E_*(X\times Y).\end{equation}
Suppose $E$ satisfies a Künneth isomorphism, i.e. (\ref{kun}) is an isomorphism for all spaces $X$ and $Y$. In other words, $E_*(-)$ is a strong symmetric monoidal functor . As per (\ref{talg}), we get an induced cartesian functor $E_*:\Top\to \kcoco$, inducing a functor
\begin{equation}\label{prehuhu}E_*:\NGrRing(\Top)\to \NGrRing(\kcoco).\end{equation} Thus for a topological ring $S$, $(E_*(B^nS))_{n\geq 0}$ is a graded $E_*$-coalgebraic ring (Definition \ref{ladef}). For $S$ discrete this was discussed by Ravenel and Wilson \cite{rw80}.\\

As a particularly simple case, take $E=Hk$, the Eilenberg-Mac Lane spectrum of the commutative ring $k$. If $k$ is a field, then $Hk$ satisfies a Künneth isomorphism, but we will not need this in what follows.

Considering sets as discrete topological spaces and modules as graded modules concentrated in degree zero, the functor $(Hk)_*:\Top\to (Hk)_*\text{-}\textbf{CoCoalg}$ restricts to the functor \begin{equation}\label{kora}k[-]:\Set \to \kkcoco\end{equation} which maps a set $X$ to the free $k$-module $k[X]$ together with the comultiplication obtained by extending linearly the diagonal map on basis elements, $\Delta(x)=x\otimes x$ for $x\in X$. In other words, the cartesian functor (\ref{kora}) is obtained from the strong symmetric monoidal free functor $k[-]:\Set \to \kmod$ by passing to cocommutative comonoids.\footnote{Remark that (\ref{kora}) cannot rightly be called a ``free functor'', since it is not the left adjoint to the ``underlying set'' functor: rather, it is the left adjoint to the ``set of group-like elements'' functor, which maps a coalgebra $C$ to the set of elements $c\in C$ such that $\Delta(c)=c\otimes c$. } 

Thus, if $S$ is a discrete ring, then $(Hk)_*(S)$ is the group ring\footnote{More accurately, to be coherent with the naming convention for this kind of object, we should say ``the ring $k$-coalgebraic ring $k[S]$''.} $k[S]$ as a $k$-coalgebraic ring concentrated in degree zero. This is the correct way of characterizing all the structure of the object $k[S]$: in particular, of characterizing the distributivity of the operations coming from the sum and multiplication of $S$.\\ 

It should be noted that there are not that many topological rings, since if $G$ is an abelian topological group then $G$ is homotopy equivalent to $\prod\limits_{n\geq 0} K(\pi_n G,n)$.

\subsection{Simplicial modules} \label{simpex}

Consider the closed symmetric monoidal category of $k$\hyp{}modules $\V=\kMod$ with the tensor product $\otimes_k$. Then 
$\CMon(\kMod)^\aug$ is the category $\kecalg$ of commutative augmented $k$-algebras, i.e. commutative $k$-algebras $A$ with a $k$-algebra homomorphism $A\to k$. The functor \begin{equation}\label{simpla}B_\bullet: \kecalg\to s(\kecalg)\end{equation} coincides with the reduced (i.e. with coefficients in the trivial $k$-module) simplicial Hochschild functor $HH_\bullet(-,k)$, as is readily verified from the definitions. However, we cannot go any further with this example: we do not have a natural choice of a functor $F$ from simplicial sets to $k$-modules. To put it differently, we have not been able to find a choice of a cosimplicial module which would yield an interesting realization of a simplicial module into a module. 

We thus shift our attention to $\V=s(\kMod)$, the closed symmetric monoidal category of simplicial $k$-modules with pointwise tensor product. It admits a natural strong symmetric monoidal functor from $s\Set$: the free simplicial $k$-module functor \[F=k[-]:s\Set \to s(\kmod).\] It has as right adjoint the functor that forgets the module structure at each level.

We claim that, as for simplicial sets (Example \ref{ej1}), the induced geometric realization functor \[|-|_{k[\Delta^\bullet]}:s^2(\kmod)\to s(\kmod)\]is the diagonal functor. The proof is very similar, but uses instead the \emph{enriched} density theorem \cite[3.72]{kelly}.

As remarked in \ref{bisimpli}, the functor $B_\bullet: s(\kecalg)\to s^2(\kecalg)$ is, degreewise, the simplicial bar construction (\ref{simpla}). After geometric realizing (taking diagonals), we obtain a functor
\[B: s(\kecalg)\to s(\kecalg)\]
which is the $k$-linear counterpart to the construction for simplicial abelian groups (\ref{bsset}). It is weakly homotopy equivalent to the algebraic $\bar{W}$-construction of Eilenberg and Mac Lane (see \cite{dtp}, A.14 for the definition of $\bar{W}$ and A.20 for a proof of $BA\simeq \bar{W}A$).

Proposition \ref{natural2} applied to $F=\id_{s\Set}$ and to $G=k[-]$ gives a natural isomorphism
\[\xymatrix{s\CMon^\aug \ar[d]_-{k[-]} \ar[r]^-{B} &  \ar[d]^-{k[-]}  s\CMon^\aug\\ s(\kecalg) \ar[r]_-{B} & s(\kecalg). \xtwocell[-1,-1]{}\omit{<0>}}\]

We cannot naively apply Theorem \ref{elteo} to this situation because $\kMod$ is not cartesian. However, we can apply (\ref{elteocoms}) to obtain a natural isomorphism
\[ \xymatrix{ s\Ring \ar[d]_-{k[-]} \ar[r]^-{B^*} & s\NGrRing \ar[d]^-{k[-]} \\ s\Ring(\kkcoco) \ar[r]_-{B^*} & s\NGrRing(\kkcoco). \xtwocell[-1,-1]{}\omit{<0>} }\]
The $B^*$ in the upper line is (\ref{bring}). Thus, if $S$ is a simplicial ring, then $k[B^*S] \cong B^*k[S]$ as simplicial graded coalgebraic rings. Let us see how the graded multiplication in $B^*k[S]$ passes to homotopy, just as we passed to homology in Example \ref{topor}. For this we need $k$ to be a field, so that the functor \[\pi_*:s(\kmod)\to \Gr(\kmod)\] is strong symmetric monoidal. Indeed, we can decompose $\pi_*$ as the normalized Moore functor $N$ (which is lax symmetric monoidal with the shuffle product; see Example \ref{degeme} for more details) followed by the homology functor. The Eilenberg-Zilber and algebraic Künneth theorems give the result. Thus, we get a functor \[\pi_*:s\NGrRing(\kkcoco) \to \NGrRing(\kkcoco).\]

Let $A\in s\Ring(\kkcoco)$ be a constant simplicial object. By neglect of structure, $A$ is an augmented commutative $k$-algebra, where the augmentation is given by the counit. Then the homotopy of $B^nA$ gives $HH_*^{[n]}(A,k)$, Pirashvili's \cite{pirashvili} higher order reduced Hochschild homology of $A$, as noted in \cite[Section 3.1]{livernet-richter}. The multiplication of $A$ induces a graded multiplication
\begin{equation}\label{grapir}HH_*^{[n]}(A,k) \otimes HH_*^{[m]}(A,k) \to HH_*^{[n+m]}(A,k).\end{equation}
Now let $S$ be a constant simplicial ring and consider $A=k[S]$. As seen in Example \ref{topor}, $B^*S$ is a model for the Eilenberg-Mac Lane graded simplicial ring $K(S,*)$. Thus, we have an isomorphism \[HH_*^{[*]}(k[S],k)\cong k[K(S,*)]\] and the graded multiplication (\ref{grapir}) corresponds to the graded multiplication (\ref{barsim}) corresponding to the cup product in cohomology. 

The reader might want to jump to Section \ref{subto} where we analyze the analogous phenomena happening in \emph{topological} Hochschild homology. 

\subsection{Differential graded modules}\label{degeme} Let $\W=\kdgm$ be the closed symmetric monoidal category of non-negatively graded chain complexes of $k$-modules with the tensor product $\otimes_k$. We would like the normalized Moore functor $N:s(\kMod)\to \kdgm$ \cite[III.2]{goerss-jardine} to play the role of our functor called $G$; however, it is not \emph{strong} symmetric monoidal, it is merely colax (with an Alexander-Whitney map) and symmetric lax (with a shuffle product map, also called ``Eilenberg-Zilber map'')\footnote{It is even \emph{bilax}: \cite[Corollary 5.7]{aguiar}.} which is why we do not take this example too far.

The functor $N$ can be seen as an extrinsic geometric realization, following Remark \ref{enrich}.\ref{ajoba}. It can be considered a sort of linearization of $|-|^e:s\Set\to \Top$. Indeed, the category of chain complexes is enriched, tensored and cotensored over modules. The cosimplicial chain complex $Nk[\Delta^\bullet]:\DDelta \to \kdgm$ begets a functor \[|-|^e_{Nk[\Delta^\bullet]}:s(\kmod)\to \kdgm.\]

The fact that $N\cong|-|^e_{Nk[\Delta^\bullet]}$ was observed in \cite[6.3]{kanfunctors} by direct computation. 
Alternatively, if we know that the Dold-Kan correspondence is an \emph{adjoint} equivalence \cite[8.4.2]{weibel}, then just as in Lemma \ref{adjunto} we get that the right adjoint to $|-|^e_{Nk[\Delta^\bullet]}$ is the functor $A\mapsto \underline{\kdgm}(Nk[\Delta^\bullet],A)$. This functor is known to be the right adjoint to $N$ in the Dold-Kan correspondence, therefore $N\cong|-|^e_{Nk[\Delta^\bullet]}$ by uniqueness of left adjoints.

The identification $N\cong |-|^e_{Nk[\Delta^\bullet]}$ is a priori a bit unsatisfactory, since to define $|-|^e_{Nk[\Delta^\bullet]}$ we use the cosimplicial chain complex $Nk[\Delta^\bullet]$ which depends on $N$. However, $Nk[\Delta^\bullet]$ can be defined independently: for example, it is the cellular chain complex on the cosimplicial space $|\Delta^\bullet|^e$. 
 
At any rate, the cosimplicial chain complex $Nk[\Delta^\bullet]:\DDelta \to \kdgm$ yields an interesting geometric realization
\[|-|_{Nk[\Delta^\bullet]}:s(\kdgm)\to \kdgm.\]
Indeed, $|-|_{Nk[\Delta^\bullet]}$ coincides with the functor $C:s(\kdgm)\to \kdgm$ which is the composition of the functors ``normalized Moore in each internal degree'', yielding a bicomplex, and the ``totalization of a bicomplex'' functor. One can prove that $C\cong|-|_{Nk[\Delta^\bullet]}$ by a computation similar to the one by Kan which proves that $N\cong |-|^e_{Nk[\Delta^\bullet]}$. 

Gugenheim and May \cite[A.2]{dtp} call the functor $C$ \emph{condensation}, and they  prove in Proposition A.3 that it is a normal lax symmetric monoidal functor (not strong), via a suitable totalized shuffle product, and colax, via a suitable Alexander-Whitney map.

Note that the induced bar construction \[B:\kecdga \to \kecdga\]
is the usual bar construction of a commutative differential graded augmented $k$-algebra, as introduced in \cite[Theorem 11.1]{onthe}; 
see also \cite[Page 69]{dtp}. Its multiplication is the ``shuffle product'', since it is induced by the lax structure of $N$. 

Proposition \ref{natural2} applies, giving a natural transformation
\[\xymatrix{s(\kecalg) \ar[d]_-N \ar[r]^-B &  \ar[d]^-N  s(\kecalg)\\ \kecdga \ar[r]_-B & \kecdga. \xtwocell[-1,-1]{}\omit{<0>}}\]
It is not a natural isomorphism since $N$ is not strong. Theorem \ref{elteo} does not apply, since $\kecdga$ is not cartesian, and Corollary \ref{elteocom} does not apply either, since $N$ is not strong.
\\

As a final remark, recall Dold and Puppe's version of the Eilenberg-Zilber theorem \cite[IV.2.4]{goerss-jardine}, which states that there is a natural quasi-isomorphism \[N(\diag X_{\bullet,\bullet}) \to CN(X_{\bullet,\bullet})\]
for a bisimplicial module $X_{\bullet, \bullet}$. Rephrasing it in the terminology adopted here, there is a natural quasi-isomorphism $N|-|\Rightarrow |N-|$. It would be interesting to check whether it is a \emph{monoidal} isomorphism, by comparison with Theorem \ref{laposta}.

\subsection{Brave new algebra} \label{bravenew}

\subsubsection{Spectra} \label{spectro}
Take $\W$ to be the closed symmetric monoidal category $\SMod$ of $S$-modules \cite{ekmm}, which we will refer to simply as ``spectra''. The monoidal structure of $\SMod$ is given by the smash product $\wedge=\wedge_\bbS$ as the tensor product, and the sphere spectrum $\bbS$ as the unit.

We take as a functor $G$ the functor \begin{equation}\label{sipa}\Sigma^\infty_+:\Top \to \SMod\end{equation} which maps a topological space $X$ to the suspension spectrum on the space $X$ with an added disjoint basepoint. It is strong symmetric monoidal \cite[II.1.2]{ekmm}, and it is a left adjoint \cite[Page 39]{ekmm}. We therefore consider $\SMod$ endowed with the cosimplicial spectrum $\Sigma^\infty_+|\Delta^\bullet|^e:\DDelta \to \SMod$ (recall from Example \ref{topor} that $|-|^e:s\Set\to \Top$ is the classical ``extrinsic'' geometric realization functor). 

It should be noted that if $X\in \Top$ and $C\in \SMod$, then $C\wedge \sip X$ defines the standard tensoring of $\SMod$ over $\Top$ \cite[III.1.1]{ekmm}. Therefore, the induced geometric realization \begin{equation}\label{geosp} |-|_{\Sigma^\infty_+|\Delta^\bullet|^e}\eqqcolon|-|:s(\SMod) \to \SMod\end{equation} in our sense coincides with the one in \cite[X.1.1]{ekmm}.

By Theorem \ref{realstrong}, $|-|$ is strong symmetric monoidal. This appears in \cite[X.1.4]{ekmm}. Theorem \ref{laposta} applies to prove that the natural isomorphism $\sip|X_\bullet|\cong |\sip X_\bullet|$, which appears in \cite[X.1.3.i]{ekmm}, is monoidal. This has not, to our knowledge, explicitly appeared in the literature.\\

We could apply Corollary \ref{elteocom} to the strong symmetric monoidal functor $G$ of (\ref{sipa}) right now. But instead, let us first delve into general $R$-modules and apply the machinery there (see (\ref{eloje})): it will give a more general result.

\subsubsection{$R$-modules and extension of scalars}Let $R$ be a commutative ring spectrum and take $\W$ to be the closed symmetric monoidal category $\RMod$ of $R$-modules, with smash product $\wedge=\wedge_R$ as tensor product, and as unit the $R$-module $R$. 

We take as a functor $G$ the strong symmetric monoidal functor of extension of scalars \[R\wedge_\bbS-:\SMod \to \RMod,\]
whose right adjoint is the restriction of scalars functor. Thus $\RMod$ is endowed with the cosimplicial $R$-module $R\wedge_\bbS \sip|\Delta^\bullet|^e$.

Furthermore, if $R'$ is another commutative ring spectrum and $f:R\to R'$ is a morphism, then in just the same fashion we obtain a strong symmetric monoidal functor $G$
\[R' \wedge_R -:\RMod\to R'\Mod\]
endowing $R'\Mod$ with the cosimplicial $R'$-module $R' \wedge_R R\wedge_\bbS \sip|\Delta^\bullet|^e \cong R'\wedge_\bbS \sip |\Delta^\bullet|^e$.

Corollary \ref{elteocom} gives a natural isomorphism comparing the iterated bar construction together with its graded structure, whenever carried out in $\RMod$ or $R'\Mod$:
\[\xymatrix{\Ring(\rcoco)\ar[d]_-{R'\wedge_R -} \ar[r]^-{B^*} &  \ar[d]^{R'\wedge_R -}  \NGrRing(\rcoco) \\ \Ring(\rrcoco) \ar[r]_-{B^*} & \NGrRing(\rrcoco). \xtwocell[-1,-1]{}\omit{<0>}  }\]

We will now analyze what the bar construction in this context actually is.

\subsubsection{Topological Hochschild homology} \label{subto} The simplicial bar construction in the category $\RMod$, \[B_\bullet: \RCAlg^\aug\to s(\RCAlg^\aug),\] coincides with $THH^R_\bullet(-,R)$, the reduced simplicial topological Hochschild homology functor, by the same verification as in Section \ref{simpex} for Hochschild homology.

Therefore,
\[B: \RCAlg^\aug\to \RCAlg^\aug\]
is an explicit model for the reduced topological Hochschild homology functor $THH^R(-,R)$, as noted in \cite[IX.2]{ekmm}. 

Note that if one wants this $THH^R(A,R)$ to coincide with the derived smash product $R \wedge_{A^e}^L A$ using the model structure in $\RMod$ set in \cite{ekmm}, it is necessary that $R$ be a q-cofibrant commutative $S$-algebra and that $A$ be a q-cofibrant commutative $R$-algebra. This is enough by a slight modification of Theorem 2.6 in \cite{ekmm}.\\

The iterations of $B$,
\begin{equation}\label{barit} B^n:\RCAlg^\aug\to \RCAlg^\aug\end{equation} for $n\geq 0$ are an explicit model for higher reduced topological Hochschild homology $THH^{R,[n]}(-,R)$ as considered e.g. in \cite{blprz}.

Indeed, $THH^R(A,R)$ can be expressed as $S^1\odot A$, where $\odot$ denotes the tensoring of the category $\recalg$ over pointed topological spaces \cite[7.1]{kuhn}, and its higher version is \[THH^{R,[n]}(A,R)= S^n \odot A.\] We obtain natural isomorphisms
\begin{align*}
B^2(A)&=THH^R(THH^R(A,R),R) = S^1\odot (S^1\odot A) \cong \\ &\cong (S^1\wedge S^1) \odot A \cong S^2\odot A = THH^{R,[2]}(A,R)
\end{align*}
and similarly for higher powers.\\

We can now apply Corollary \ref{elteocom} to the functor 
\begin{equation}\label{eloje}G=R[-]=R\wedge_\bbS \sip:\Top\to R\Mod,\end{equation}and obtain a natural isomorphism
\[\xymatrix{ \Ring(\Top) \ar[d]_-{R[-]} \ar[r]^-{B^*} & \NGrRing(\Top) \ar[d]^-{R[-]} \\ \Ring(\rcoco) \ar[r]_-{B^*} & \NGrRing(\rcoco). \xtwocell[-1,-1]{}\omit{<0>} }\]

We thus get, for $A\in \Ring(\rcoco)$, a graded multiplication in higher $THH$:
\begin{equation}\label{elprod}THH^{R,[n]}(A,R) \wedge_R THH^{R,[m]}(A,R)\to THH^{R,[n+m]}(A,R)\end{equation}
and if $A=R[S]$ for $S\in \Ring(\Top)$, then we get a natural isomorphism
\begin{equation}\label{loultimo}THH^{R,[*]}(R[S],R) \cong R[B^* S]\end{equation}
of graded ring objects of $\rcoco$, i.e. of coalgebraic rings in $\RMod$ (Definition \ref{ladef}). As noted in Example \ref{topor}, when $S$ is discrete the graded topological ring $B^*S$ is a model for the Eilenberg-Mac Lane spaces $K(S,*)$.

Observe that by means of Remark \ref{nojodas}, we get for every $n\geq 0$ an isomorphism
\[THH^{R,[n]}(R[G], R)\cong R[B^nG]\]
of bicommutative $R$-Hopf algebras, natural in the topological abelian group $G$. For $n=1$ and $R=\bbS$ there is a similar well-known unreduced result: $THH(\bbS[G]) \cong \bbS[B^{cy}G]$, 
where $B^{cy}$ stands for the cyclic bar construction.\\ 

Let us pass to homology. Let $E\in \RCAlg$ be a field, meaning that $E_*=\pi_*E$ is a graded field, i.e. every graded module over it is free. Then the $E$-homology functor $E_*:\RMod\to E_*\Mod$ is strong symmetric monoidal, since the field hypothesis guarantees a Künneth isomorphism. It induces a functor
\begin{equation}\label{huhu}E_*:\NGrRing(\rcoco) \to \NGrRing(\kcoco).\end{equation}
Thus, $E_*(THH^{R,[*]}(R[S],R))$ 
is an $E_*$-coalgebraic ring. If in particular we take $R$ to be $\bbS$, we can recover the coalgebraic ring of Ravenel-Wilson (see Example \ref{topor}) in a different guise. Indeed, precomposing (\ref{huhu}) with \[\bbS[-]:\NGrRing(\Top)\to \NGrRing(\rcoco)\] gives the $E$-homology of topological rings (\ref{prehuhu}). Then, if $S\in\Ring(\Top)$, we get an isomorphism of $E_*$-coalgebraic rings \begin{equation}E_*(THH^{\bbS,[*]}(\bbS[S],\bbS)) \cong E_*(B^*S).\end{equation} 

\appendix
\section{Monoidal category theory} 

We recall some results we will use from monoidal category theory. A modern, comprehensive reference is \cite{aguiar}. 

\subsection{General facts} \label{general}

\bi
\item In monoidal categories $\V$ we can consider \emph{(co)monoids} and morphisms of these. Thus we get categories $\Mon(\V)$ and $\Comon(\V)$ which are in turn symmetric monoidal provided $\V$ is symmetric. The symmetry is needed to have a natural multiplication on the tensor product of monoids. We will henceforth focus on symmetric monoidal categories, since this is the setting we need in the body of the paper.

In symmetric monoidal categories, we can also consider (co)commutative \linebreak (co)monoids, and thus get symmetric monoidal categories $\CMon(\V)$ and \linebreak $\CoComon(\V)$.

The unit for these categories of monoids and comonoids is the unit of $\V$, which has unique commutative monoid and cocommutative comonoid structures. It is an initial object for $\Mon(\V)$ and $\CMon(\V)$, and a final object for $\Comon(\V)$ and for $\CoComon(\V)$. 
\item
A \emph{lax monoidal functor} $F:\V\to \W$ between monoidal categories is the data of a functor $F$ plus morphisms \[\nabla_{V,V'}:F(V)\otimes F(V')\to F(V\otimes V')\] 
natural in $V,V'\in \V$, and a unit map \[\nabla_0: \1_\W\to F(\1_\V)\] satisfying associativity and unitality conditions. A lax monoidal functor is \emph{normal} if $\nabla_0$ is an isomorphism (see Section \ref{moncaug} for an application of this concept). If the categories are symmetric monoidal, then a lax monoidal functor is \emph{symmetric} if it commutes with the respective symmetries. 

Dually, a \emph{colax monoidal functor} $F:\V\to \W$ between monoidal categories is the data of a functor $F$ plus natural morphisms \[\Delta_{V,V'}:F(V\otimes V')\to F(V)\otimes F(V')\] and a counit map \[\Delta_0: F(\1_\V) \to \1_\W\] satisfying coassociativity and counitality conditions. 
If the categories are symmetric monoidal, then a colax monoidal functor is \emph{symmetric} if it commutes with the respective symmetries.
\item \label{induces} An important feature of (co)lax monoidal functors is that they send (co)monoids to (co)monoids. More precisely, if $F:\V\to \W$ is a (co)lax monoidal functor between symmetric monoidal categories, we get an induced monoidal functor $F:\Mon(\V)\to \Mon(\W)$ (resp. $F:\Comon(\V)\to \Comon(\W)$). For example, if $M\in \Mon(\V)$ has multiplication $\mu$ and $F$ is lax monoidal, then the multiplication in $F(M)$ is given by the composition
\[\xymatrix{F(M)\otimes F(M)\ar[r]^-{\nabla_{M,M}} & F(M\otimes M) \ar[r]^-{F(\mu)} & F(M).}\]
If moreover $F$ is symmetric (co)lax, we get an induced symmetric monoidal functor $F:\CMon(\V)\to \CMon(\W)$ (resp. $F:\CoComon(\V)\to \CoComon(\W)$). 
\item A \emph{strong monoidal functor} $F:\V\to \W$ between monoidal categories is a lax monoidal functor such that $\nabla_{V,V'}$ and $\nabla_0$ are isomorphisms (in particular, it is normal). Thus a strong monoidal functor is also a colax monoidal functor, by inverting the structure morphisms. Therefore, if the categories are symmetric then $F$ induces monoidal functors $F:\Mon(\V)\to \Mon(\W)$ and $F:\Comon(\V)\to \Comon(\W)$. If $F$ is symmetric, then it also induces symmetric monoidal functors between the categories of (co)commutative (co)monoids.
\item \emph{The Eckmann-Hilton argument.} Let $\V$ be a symmetric monoidal category. There is an equivalence of symmetric monoidal categories (i.e. a strong symmetric monoidal functor which is an equivalence of categories)\[\Mon(\Mon(\V))\cong\CMon(\V).\]
The symmetric monoidal category $\CMon(\V)$ also coincides with the symmetric monoidal categories $\Mon(\CMon(\V))$, $\CMon(\Mon(\V))$, $\CMon(\CMon(\V))$.
\ei

\subsubsection{2-categorical remarks} \label{2cats}
\bi
\item One can define natural transformations of (symmetric) lax, colax, strong monoidal functors \cite[3.1.2]{aguiar}. For each of these alternatives there is a corresponding 2-category \cite[3.3.3, C.1.1]{aguiar}. The natural transformations of (co)lax monoidal functors are called \emph{monoidal transformations}. The natural transformations of (symmetric) strong monoidal functors and normal (co)lax monoidal functors are the natural transformations of the underlying (co)lax monoidal functors.

Let us make the lax case explicit. A monoidal transformation $\xymatrix@C+1pc{\V \rtwocell<4>^{F}_G{\;\tau}  & \W}$ of lax monoidal functors is a natural transformation such that the following diagrams commute for every $V,V'\in \V$.
\[\xymatrix{FV \otimes FV' \ar[r]^-{\nabla_{V,V'}^F} \ar[d]_-{\tau_V\otimes \tau_{V'}} & F(V\otimes V') \ar[d]^-{\tau_{V\otimes V'}} \\ GV \otimes GV' \ar[r]_-{\nabla_{V,V'}^G} & G(V\otimes V').} \hspace{1cm}
\xymatrix{ & F\1_\V \ar[dd]^-{\tau_{\1_\V}} \\ \1_\W \ar[ru]^-{\nabla_0^F} \ar[rd]_-{\nabla_0^G} \\ & G\1_\V} \] 
\item Let $\xymatrix@C+1pc{\V \rtwocell<4>^{F}_G{\;\tau}  & \W}$ be a 2-cell in the 2-category of symmetric monoidal categories, lax symmetric monoidal functors and monoidal transformations. There is an induced 2-cell \begin{equation}\label{2cosa} \xymatrix@C+1pc{\CMon(\V) \rtwocell<6>^{F}_G{\;\tau}  & \CMon(\W)}\end{equation} in the same 2-category. The fact that $\tau$ is a monoidal transformation and not merely a natural transformation guarantees that if $A\in \CMon(\V)$, then $\tau_A:FA\to GA$ is a morphism of commutative monoids in $\W$.
\ei

\subsection{Cartesian monoidal categories} \label{cartesiansect}

\bi
\item A symmetric monoidal category $\C$ is \emph{cartesian} if the monoidal structure defines category-theoretical binary products: the tensor product is a categorical product and the unit is a final object.

Any category where finite products exist is cartesian monoidal, after a definite, functorial choice of finite products.

In a cartesian category, more structure is available than in a mere symmetric monoidal category: we have projections and diagonal maps.
\item \label{strcar} A strong symmetric monoidal functor between cartesian categories is called \emph{cartesian}. They are also called ``finite product-preserving functors'', because they send the chosen product diagram $A\leftarrow A\times B \to B$ to a product diagram. 
\item \label{prop13} If $\V$ is a symmetric monoidal category, then the symmetric monoidal category 
$\CoComon(\V)$ is cartesian.
\item \label{prop12} Conversely, if $\C$ is a cartesian category, then every object of $\C$ is uniquely a cocommutative comonoid for the product: $\C=\CoComon(\C)$, thanks to the existence of diagonal maps $C\to C\times C$. 

In particular, if we are given a colax symmetric monoidal functor $\C\to \V$ from a cartesian category to a symmetric monoidal category, then we get an induced colax symmetric monoidal functor \begin{equation}\label{laec} \C\to \CoComon(\V)\end{equation}
which is cartesian if $\C\to \V$ is strong.
\item \label{talg} In a cartesian category $\C$, one can internalize algebras for a Lawvere theory (also called ``algebraic theories'' or ``finite product theories''; see \cite[Chapter 3]{borceux2}), and these categories of algebras are in turn cartesian. A cartesian functor between cartesian categories induces cartesian functors between the categories of algebras. 

We will be considering abelian groups objects $\Ab(\C)$, ring objects $\Ring(\C)$,  and $\N$-graded ring objects $\NGrRing(\C)$.
\ei
\subsection{Hopf monoids} \label{hopfmo}
\bi
\item A \emph{Hopf monoid} \cite[1.2.5]{aguiar} in a symmetric monoidal category $\V$ is an object of $\V$ together with compatible monoid and comonoid structures, plus an antipode morphism: the prototypical example is a Hopf algebra, which is a Hopf monoid in a category of modules over a commutative ring. A Hopf monoid is \emph{bicommutative} if it is commutative and cocommutative; these form a symmetric monoidal category denoted by $\biHopf(\V)$. It turns out that in this case, \begin{equation}\label{biho}\biHopf(\V) \cong \Ab(\CoComon(\V)),\end{equation} the category of abelian group objects in the cartesian category $\CoComon(\V)$ (\cite[20.3.4]{mayponto}, see also \cite[1.22.iii]{aguiar}). In particular, $\biHopf(\V)$ is cartesian. 
\item A strong symmetric monoidal functor $F:\V\to \W$ between symmetric monoidal categories induces a cartesian functor $F:\biHopf(\V)\to \biHopf(\W)$. 
Indeed, we can break it up in two parts: first, $F$ induces a cartesian functor between the categories of cocommutative monoids. Then, this functor induces a functor in the categories of abelian group objects. We are done by the characterization (\ref{biho}).
\item There is a forgetful functor $\biHopf(\V)\to \CMon(\V)^\aug$, where the augmentation to be considered is the counit of the given Hopf monoid. This forgetting procedure is compatible with strong symmetric monoidal functors $\V\to \W$.
\ei

\subsection{Overcategories} \label{augmented}

\subsubsection{Augmented objects} \label{augmented2}

\bi
\item Let $\C$ be a category with an initial object $I$. Define the category $\C^\aug$ of \emph{augmented objects} of $\C$ to be the slice category $\C^{/I}$, i.e. the category of objects of $\C$ with a chosen morphism to $I$, with the obvious morphisms.

\item In $\C^\aug$, the object $\xymatrix{I\ar[r]^-\id &I}$ is a zero object  (i.e. both initial and final). The forgetful functor $\C^\aug \to \C$ is an equivalence of categories if and only if $\C$ has a zero object. In particular, $(\C^\aug)^\aug\cong \C^\aug$.
\ei

\subsubsection{Augmented objects in categories of monoids and abelian groups} \label{moncaug}
\bi
\item As we have already observed, $\Mon(\V)$ and $\CMon(\V)$ have $\1$ as their initial object. Thus, we can consider $\Mon(\V)^\aug$ and $\CMon(\V)^\aug$, the categories of \emph{augmented (commutative) monoids}. They have induced symmetric monoidal structures.

\item Let $F:\V\to \W$ be a normal lax symmetric monoidal functor. There is an induced normal lax symmetric monoidal functor $F:\CMon(\V)^\aug\to \CMon(\W)^\aug$. Normality is needed so that an augmentation of a monoid $A\to \1$ in $\V$ gets mapped to an augmentation $FA \to F\1 \stackrel{\cong}{\leftarrow} \1$ of $FA$ in $\W$. 

If $\tau$ is a monoidal transformation between normal lax symmetric monoidal functors, there is an induced 2-cell like (\ref{2cosa}) but on augmented commutative monoids:
\begin{equation}\label{cosados} \xymatrix@C+1pc{\CMon(\V)^\aug \rtwocell<6>^{F}_G{\;\tau}  & \CMon(\W)^\aug}\end{equation}
Similar remarks hold for augmented monoids, normal lax monoidal functors and monoidal transformations.
\item If $\C$ is a cartesian category then $\1$ is a final object, and $\Ab(\C)$ has $\1$ as a zero object. Thus, $\Ab(\C)^\aug\cong \Ab(\C)$.
\ei

\subsection{Simplicial objects} \label{simplicial}

\bi
\item If $\C$ is any category, we denote by $s\C$ its category of simplicial objects, i.e. the functor category $\Fun(\DDelta^\op,\C)$ where $\DDelta$ is the category of non-empty finite ordinals and order-preserving maps. If $\C$ is (symmetric) monoidal or cartesian, then so is $s\C$, with pointwise tensor product: \[(X_\bullet \otimes Y_\bullet)_n=X_n\otimes Y_n.\] Its unit is the constant simplicial object $\1$. Any functor $F:\C\to \D$ induces a functor $F:s\C\to s\D$, which is cartesian or (symmetric, normal) lax, colax, strong monoidal if the original one is.

Let $s^2\C=\Fun(\DDelta^\op \times \DDelta^\op,\C)$ be the category of bisimplicial objects in $\C$. By adjunction, we can identify $s^2\C$ with $s(s\C)$.  
Thus, we will think of a bisimplicial object $X_{\bullet, \bullet}\in s^2\C$ as a simplicial object $([n]\mapsto X_{n,\bullet})$ in $s\C$.

\item For $\V$ a symmetric monoidal category, there is an equivalence of categories \linebreak $s(\V^{/\1})\cong (s\V)^{/\1}$ (the slice category notation was explained in \ref{augmented2}). Similarly, for $\C$ a category with an initial object $I$, we have $s(\C^\aug)\cong (s\C)^\aug$.

\item Let $T$ be the cartesian category associated to a Lawvere theory (cf. the last bullet point of \ref{cartesiansect}) and $\C$ be a cartesian category. Then the category of $\C$-valued $T$-algebras is the category of cartesian functors $T\to\C$.

By adjunction, 
the category of $s\C$-valued $T$-algebras coincides with the category of functors from $\DDelta^\op$ to $\C$-valued $T$-algebras. In symbols,
\[\mathrm{CFun}(T,s\C)=\Fun(\DDelta^\op,\mathrm{CFun}(T,\C))\]
where $\mathrm{CFun}$ denotes cartesian functors.

In particular, $s(\Ab(\C))=\Ab(s\C)$, and similarly for rings and graded rings. 

\item This also works for ((co)commutative) (co)monoids in a (symmetric) monoidal category $\V$. 
Thus, $s(\Mon(\V))=\Mon(s\V)$ and similarly for commutative monoids and (co)commutative comonoids.
\ei

\bibliographystyle{alpha}
\bibliography{main}

\newcommand{\etalchar}[1]{$^{#1}$}
\begin{thebibliography}{EKMM97}

\bibitem[AM10]{aguiar}
Marcelo Aguiar and Swapneel Mahajan.
\newblock {\em Monoidal functors, species and {H}opf algebras}, volume~29 of
  {\em CRM Monograph Series}.
\newblock American Mathematical Society, Providence, RI, 2010.
\newblock With forewords by Kenneth Brown and Stephen Chase and Andr{\'e}
  Joyal.

\bibitem[BLP{\etalchar{+}}15]{blprz}
Irina Bobkova, Ayelet Lindenstrauss, Kate Poirier, Birgit Richter, and Inna
  Zakharevich.
\newblock On the higher topological {H}ochschild homology of {$\Bbb F_p$} and
  commutative {$\Bbb F_p$}-group algebras.
\newblock In {\em Women in topology: collaborations in homotopy theory}, volume
  641 of {\em Contemp. Math.}, pages 97--122. Amer. Math. Soc., Providence, RI,
  2015.

\bibitem[Bor94]{borceux2}
Francis Borceux.
\newblock {\em Handbook of categorical algebra. 2}, volume~51 of {\em
  Encyclopedia of Mathematics and its Applications}.
\newblock Cambridge University Press, Cambridge, 1994.
\newblock Categories and structures.

\bibitem[EKMM97]{ekmm}
A.~D. Elmendorf, I.~Kriz, M.~A. Mandell, and J.~P. May.
\newblock {\em Rings, modules, and algebras in stable homotopy theory},
  volume~47 of {\em Mathematical Surveys and Monographs}.
\newblock American Mathematical Society, Providence, RI, 1997.
\newblock With an appendix by M. Cole.

\bibitem[EML53]{onthe}
Samuel Eilenberg and Saunders Mac~Lane.
\newblock On the groups of {$H(\Pi,n)$}. {I}.
\newblock {\em Ann. of Math. (2)}, 58:55--106, 1953.

\bibitem[GJ99]{goerss-jardine}
Paul~G. Goerss and John~F. Jardine.
\newblock {\em Simplicial homotopy theory}, volume 174 of {\em Progress in
  Mathematics}.
\newblock Birkh\"auser Verlag, Basel, 1999.

\bibitem[GM74]{dtp}
V.~K. A.~M. Gugenheim and J.~Peter May.
\newblock {\em On the theory and applications of differential torsion
  products}.
\newblock American Mathematical Society, Providence, R.I., 1974.
\newblock Memoirs of the American Mathematical Society, No. 142.

\bibitem[Kan58]{kanfunctors}
Daniel~M. Kan.
\newblock Functors involving c.s.s. complexes.
\newblock {\em Trans. Amer. Math. Soc.}, 87:330--346, 1958.

\bibitem[Kel05]{kelly}
G.~M. Kelly.
\newblock Basic concepts of enriched category theory.
\newblock {\em Repr. Theory Appl. Categ.}, (10):vi+137, 2005.
\newblock Reprint of the 1982 original [Cambridge Univ. Press, Cambridge;
  MR0651714].

\bibitem[Kuh04]{kuhn}
Nicholas~J. Kuhn.
\newblock The {M}c{C}ord model for the tensor product of a space and a
  commutative ring spectrum.
\newblock In {\em Categorical decomposition techniques in algebraic topology
  ({I}sle of {S}kye, 2001)}, volume 215 of {\em Progr. Math.}, pages 213--236.
  Birkh\"auser, Basel, 2004.

\bibitem[Lac00]{lack}
Stephen Lack.
\newblock A coherent approach to pseudomonads.
\newblock {\em Adv. Math.}, 152(2):179--202, 2000.

\bibitem[LR11]{livernet-richter}
Muriel Livernet and Birgit Richter.
\newblock An interpretation of {$E_n$}-homology as functor homology.
\newblock {\em Math. Z.}, 269(1-2):193--219, 2011.

\bibitem[May72]{gils}
J.~P. May.
\newblock {\em The geometry of iterated loop spaces}.
\newblock Springer-Verlag, Berlin-New York, 1972.
\newblock Lectures Notes in Mathematics, Vol. 271.

\bibitem[Mil67]{milgram}
R.~James Milgram.
\newblock The bar construction and abelian {$H$}-spaces.
\newblock {\em Illinois J. Math.}, 11:242--250, 1967.

\bibitem[ML70]{maclanemilgram}
Saunders Mac~Lane.
\newblock {\em The {S}teenrod Algebra and Its Applications: A Conference to
  Celebrate {N}.{E}. {S}teenrod's Sixtieth Birthday: Proceedings of the
  Conference held at the {B}attelle {M}emorial {I}nstitute, {C}olumbus, {O}hio
  {M}arch 30th—{A}pril 4th, 1970}, chapter The Milgram bar construction as a
  tensor product of functors, pages 135--152.
\newblock Springer Berlin Heidelberg, Berlin, Heidelberg, 1970.

\bibitem[ML98]{maclane}
Saunders Mac~Lane.
\newblock {\em Categories for the working mathematician}, volume~5 of {\em
  Graduate Texts in Mathematics}.
\newblock Springer-Verlag, New York, second edition, 1998.

\bibitem[MP12]{mayponto}
J.~P. May and K.~Ponto.
\newblock {\em More concise algebraic topology}.
\newblock Chicago Lectures in Mathematics. University of Chicago Press,
  Chicago, IL, 2012.
\newblock Localization, completion, and model categories.

\bibitem[Pir00]{pirashvili}
Teimuraz Pirashvili.
\newblock Hodge decomposition for higher order {H}ochschild homology.
\newblock {\em Ann. Sci. \'Ecole Norm. Sup. (4)}, 33(2):151--179, 2000.

\bibitem[Rie14]{riehl}
Emily Riehl.
\newblock {\em Categorical homotopy theory}, volume~24 of {\em New Mathematical
  Monographs}.
\newblock Cambridge University Press, Cambridge, 2014.

\bibitem[RW80]{rw80}
Douglas~C. Ravenel and W.~Stephen Wilson.
\newblock The {M}orava {$K$}-theories of {E}ilenberg-{M}ac {L}ane spaces and
  the {C}onner-{F}loyd conjecture.
\newblock {\em Amer. J. Math.}, 102(4):691--748, 1980.

\bibitem[RW77]{rw77}
Douglas~C. Ravenel and W.~Stephen Wilson.
\newblock The {H}opf ring for complex cobordism.
\newblock {\em J. Pure Appl. Algebra}, 9(3):241--280, 1976/77.

\bibitem[Shu]{shulmanmo}
Mike Shulman.
\newblock Coherence theorem for symmetric lax monoidal functors.
\newblock MathOverflow.
\newblock URL:http://mathoverflow.net/q/243442 (version: 2016-07-12).

\bibitem[Ste12]{stevenson}
Danny Stevenson.
\newblock D\'ecalage and {K}an's simplicial loop group functor.
\newblock {\em Theory Appl. Categ.}, 26:768--787, 2012.

\bibitem[Wei94]{weibel}
Charles~A. Weibel.
\newblock {\em An introduction to homological algebra}, volume~38 of {\em
  Cambridge Studies in Advanced Mathematics}.
\newblock Cambridge University Press, Cambridge, 1994.

\bibitem[Zis15]{zisman}
Michel Zisman.
\newblock Comparaison de deux diagonales pour un ensemble bisimplicial.
\newblock {\em J. Homotopy Relat. Struct.}, 10(4):1013--1021, 2015.

\end{thebibliography}

\end{document}